\documentclass[final]{siamltex}
\usepackage{color}
\usepackage{array}
\usepackage{graphicx}
\usepackage{amsmath,amssymb}
\usepackage{algorithm}
\usepackage{algpseudocode}
\usepackage{hyperref}
\usepackage{varwidth}
\definecolor{hotpink}{rgb}{0.9,0,0.5}
\usepackage{epstopdf}
\usepackage{multirow}
\usepackage{capt-of}

\numberwithin{theorem}{section}

\allowdisplaybreaks

\title{The Adaptive $s$-step Conjugate Gradient Method}

\author{
  Erin Carson\thanks{Courant Institute of Mathematical Sciences, New York University, New York, NY. 
    (\url{erinc@cims.nyu.edu}, \url{http://math.nyu.edu/\~ erinc/}).}
}

\def\uyh{\underline{\hat{\mathcal{Y}}}}
\def\yh{\widehat{\mathcal{Y}}}
\def\y{\mathcal{Y}}
\def\xhat{\widehat{x}}
\def\rhat{\widehat{r}}
\def\qhat{\widehat{q}}
\def\B{\mathcal{B}}
\def\Gh{\widehat{G}}
\def\veps{\varepsilon}
\def\norm#1{\|#1\|}

\begin{document}
\maketitle

\begin{abstract}
On modern large-scale parallel computers, the performance of Krylov subspace iterative methods is limited 
by global synchronization. This has inspired the development of $s$-step (or communication-avoiding) Krylov subspace method variants, in which iterations are computed in blocks of $s$. This reformulation can reduce the number of global synchronizations per iteration by a factor of $O(s)$, and has been shown to produce speedups in practical settings. 

Although the $s$-step variants are mathematically equivalent to their classical counterparts, they can behave quite differently in finite precision depending on the parameter $s$. If $s$ is chosen too large, the $s$-step method can suffer a convergence delay and a decrease in attainable accuracy relative to the 
classical method. This makes it difficult for a potential user of such methods - the $s$ value that minimizes the time per iteration may not be the best $s$ for minimizing the overall time-to-solution, and further may cause an unacceptable decrease in accuracy. 

Towards improving the reliability and usability of $s$-step Krylov subspace methods, in this work we derive the \emph{adaptive $s$-step CG method}, a variable $s$-step CG method where in block $k$, 
the parameter $s_k$ is determined automatically such that a user-specified accuracy is attainable. 
The method for determining $s_k$ is based on a bound on growth of the residual gap within block $k$, from which we derive a constraint on the condition numbers of the computed $O(s_k)$-dimensional Krylov subspace bases. The computations required for determining the block size $s_k$ can be performed without increasing the number of global synchronizations per block. 
Our numerical experiments demonstrate that the adaptive $s$-step CG method is able to attain up to the same accuracy as classical CG while still significantly reducing the total number of global synchronizations. 

%
\end{abstract}

%

\section{Introduction}\label{sec:introduction}

In this paper, we consider the use of Krylov subspace methods for solving large, sparse linear systems $Ax=b$, where $A\in \mathbb{R}^{n\times n}$. We will focus on the conjugate gradient (CG) method~\cite{hestenes1952methods}, which is used when $A$ is symmetric positive definite. 
Given an initial approximate solution $x_1$ and corresponding residual $r_1=b-Ax_1$, the CG method iteratively updates the approximate solution using a sequence of nested Krylov subspaces $\mathcal{K}_{1}(A,r_1)\subset \mathcal{K}_{2}(A,r_1) \subset \cdots$, where  
\[
\mathcal{K}_{i}(A,r_1) = \text{span}(r_1,Ar_1,\ldots,A^{i-1}r_1)
\]
denotes the $i$-dimension Krylov subspace with matrix $A$ and starting vector $v$. In iteration $i$, the updated approximate solution 
$x_{i+1}\in x_1 + \mathcal{K}_{i}(A,r_1)$ is chosen by imposing the Galerkin condition $r_{i+1}=b-Ax_{i+1} \perp \mathcal{K}_{i}$. 

Thus each iteration of CG requires a matrix-vector product with $A$ in order to expand the dimension of the Krylov subspace and a number of inner products to perform the orthogonalization. 
On modern computer architectures, the speed with which the matrix-vector products and inner products can be computed is limited by communication (i.e., the movement of data). This limits the potential speed
of individual iterations
attainable by an implementation of CG. 
To perform a sparse matrix-vector product in parallel, each processor must communicate entries of the source vector and/or the destination vector that it owns to neighboring processors. 
Inner products require a global synchronization, i.e., the computation can not proceed until all processors
have finished their local computation and communicated the result to other processors.
For large-scale sparse problems on large-scale machines, the cost of
synchronization between parallel processors can dominate the run-time
(see, e.g., the exascale computing report~\cite[pp.~28]{dongarra11}).

Research efforts toward removing the performance bottlenecks caused by communication in CG and other Krylov subspace methods have produced various approaches. 
One such approach are the $s$-step Krylov subspace methods (also called ``communication-avoiding'' Krylov subspace methods); for a thorough treatment of background, related work, and performance experiments, see, e.g., the theses~\cite{carson2015communication, hoemmen2010communication}. In $s$-step Krylov subspace methods, instead of performing one iteration at a time, the iterations are performed in blocks of $s$; i.e., in each iteration, the Krylov subspace is expanded by $O(s)$ dimensions by computing $O(s)$ new basis vectors and then all inner products between the new basis vectors needed for the next $s$ iterations are computed in one block operation. In this way, computing the inner products for $s$ iterations only requires a single global synchronization, decreasing the synchronization cost per iteration by a factor of $O(s)$. This approach has been shown to lead to significant speedups for a number of problems and real-world applications (see, e.g.,~\cite{mohiyuddin2009minimizing,williams2014ipdps}). 
In the remainder of the paper, we will refer to the matrices whose columns consist of the $O(s)$ basis vectors computed in each block as \emph{$s$-step basis matrices}. 
Further details of the $s$-step CG method are discussed in section~\ref{sec:sksms}. We emphasize that our use of the overloaded term ``$s$-step methods'' here differs from other works, e.g.,~\cite{cullum1974block,karush1951iterative} and~\cite[\S 9.2.7]{golub1996matrix}, in which `$s$-step method' refers to a type of restarted Lanczos procedure. 

In exact arithmetic the $s$-step CG method produces the exact same iterates as the classical CG method, but 
their behavior can differ significantly in finite precision.
In both $s$-step and classical Krylov subspace methods, rounding errors due to finite precision arithmetic have two basic effects: a decrease in 
attainable accuracy and a delay of convergence. 
It has long been known that for $s$-step Krylov subspace methods, as $s$ is increased (and so the condition numbers of the $s$-step basis matrices increase), the attainable accuracy decreases and the convergence delay increases relative to the classical CG method (see, e.g.,~\cite{chronopoulos1989s}). At the extreme, if the parameter $s$ is chosen to be too large, the $O(s)$-dimensional bases computed for each block can be numerically rank deficient and the $s$-step method can fail to converge. 

This sensitive numerical behavior poses a practical obstacle to optimizing the performance of $s$-step methods, and diminishes their usability and reliability. In a setting where the performance of CG is communication-bound, we expect that up to some point, increasing $s$ will decrease the time per iteration. If we pick $s$ only based on minimizing the time per iteration, however, we can run into problems. 
First, the finite precision error may cause a large convergence delay, negating any potential performance gain with respect to the overall runtime. Since the number of iterations required for convergence for a given $s$ value is not known a priori, choosing the $s$ value that results in the fastest time-to-solution is a difficult problem. Second, the chosen $s$ parameter may cause $s$-step CG to fail to converge to the user-specified accuracy. In this case, the particular problem is \emph{unsolvable} by the $s$-step CG method. 

Requiring the user to choose the parameter $s$ thus diminishes the practical usability of $s$-step Krylov subspace methods. 
It is therefore imperative that we develop a better understanding of the convergence rate and accuracy in finite precision $s$-step CG and other $s$-step Krylov subspace methods. Our hope is that by studying the theoretical properties of methods designed for large-scale computations in finite precision, we can develop methods and techniques that are efficient, capable of meeting application-dependent accuracy constraints, and which do not require that the user have extensive knowledge of numerical linear algebra.

Toward this goal, in this paper we develop a variable $s$-step CG method in which the $s$ values can vary between blocks and are chosen automatically by the algorithm such that a user-specified accuracy is attained. We call this the \emph{adaptive $s$-step CG method}.
The condition for determining $s$ in each block is derived using a bound on the size of the gap between the true residuals $b-Ax_i$ and the updated residuals $r_i$ computed in finite precision. 
We show that in order to prevent a loss of attainable accuracy, in a given block of iterations, the condition number of the $s$-step basis matrix should be no larger than a quantity related to the inverse of the 2-norms of the residuals computed within that block. In other words, as the method converges, more ill-conditioned bases (i.e., larger $s$ values) can be used without a loss of attainable accuracy. 
Furthermore, we show that the adaptive $s$-step CG method can be implemented in a way that does not increase the synchronization cost per block versus (fixed) $s$-step CG. 

Whereas the fixed $s$-step CG method may not attain the same accuracy as classical CG, 
our numerical experiments on a set of example matrices demonstrate that our adaptive $s$-step CG method can achieve the same accuracy as classical CG,  while still reducing the number of synchronizations required by up to over a factor of $4$. 

In section~\ref{sec:sksms} we give a brief background on $s$-step Krylov subspace methods and detail the $s$-step CG method, and in section~\ref{sec:relatedwork} we discuss related work. Section~\ref{sec:varscg} presents our main contribution; in section~\ref{sec:theory}, we derive a theoretical bound on how large the $s$-step basis matrix condition number can be in terms of the norms of the updated residuals without affecting the attainable accuracy, and in section~\ref{sec:algorithm}, we detail the adaptive $s$-step CG algorithm that makes use of this bound in determining the number of iterations to perform in each block. Numerical experiments for a small set of test matrices are presented in section~\ref{sec:experiments}. We discuss possible extensions and future work and briefly conclude in section~\ref{sec:conclusion}.

\section{Background: $s$-step Krylov subspace methods}
\label{sec:sksms}
				
The basic idea behind $s$-step Krylov subspace methods is to grow the underlying Krylov subspace $O(s)$ dimensions at a time and perform the necessary orthogonalization with a single global synchronization. Within a block of $s$ iterations, the vector updates are performed implicitly by updating the coordinates of the iteration vectors in the $O(s)$-dimensional basis spanned by the columns of the $s$-step basis matrix. The algorithmic details vary depending on the particular Krylov subspace method considered, but the general approach is the same. 
There is a wealth of literature related to $s$-step Krylov subspace methods. We give a brief summary of early related work below; a thorough bibliography can be found in~\cite[Table 1.1]{hoemmen2010communication}.

The $s$-step approach was first used in the context of Krylov subspace methods by Van Rosendale~\cite{rosendale1983minimizing}, who developed a variant of the parallel conjugate gradient method which minimizes inner product data dependencies with the goal of exposing more parallelism. 
The term ``$s$-step Krylov subspace methods'' was first used by Chronopoulos and Gear, who developed an $s$-step CG method~\cite{chronopoulos1989efficient, chronopoulos1989s}. This work was followed in subsequent years by the development of other $s$-step variants of Orthomin and GMRES~\cite{chronopoulos1990s}, Arnoldi and Symmetric Lanczos~\cite{kim1991class, kim1991efficient}, MINRES, GCR, and Orthomin~\cite{chronopoulos1991s,chronopoulos1996parallel}, Nonsymmetric Lanczos~\cite{kim1992efficient}, and Orthodir~\cite{chronopoulos2001odir}.
Walker similarly used $s$-step bases, but with the goal of improving stability in GMRES by replacing the modified Gram-Schmidt orthogonalization process with Householder QR \cite{walker1988implementation}. 

Many of the earliest $s$-step Krylov subspace methods used monomial bases (i.e., $[v, Av, A^2v,\ldots]$) for the computed $s$-step basis matrices, but it was found that convergence often could not be guaranteed for $s>5$ (see, e.g.,~\cite{chronopoulos1989s}). This motivated research into the use of other more well-conditioned bases for the Krylov subspaces, including scaled monomial bases~\cite{hindmarsh1986note}, Chebyshev bases~\cite{joubert1992parallelizable,de1991parallel,de1995reducing}, and Newton bases~\cite{bai1994newton,erhel1995parallel}.
%
%
The growing cost of communication in large-scale sparse problems has created a recent resurgence of interest in the implementation, optimization, and development of $s$-step Krylov subspace methods; see, e.g., the recent works~\cite{demmel2007avoiding,hoemmen2010communication,mohiyuddin2009minimizing,yamazaki2012communication,grigori2013ILU0,mehridehnavi2013communication,yamazaki2014improving,boman2014domain,williams2014ipdps}.

\subsection{The $s$-step conjugate gradient method}
In this work, we focus on the $s$-step conjugate gradient method (Algorithm~\ref{alg:scg}), which is equivalent to the $s$-step CG method (Algorithm~\ref{alg:scg}) in exact arithmetic. 
In the remainder of this section, we give a brief explanation of the $s$-step CG method; further details on the derivation, implementation, and performance of the $s$-step CG method can be found in, e.g.,~\cite{carson2015communication, hoemmen2010communication}. 
				
\begin{algorithm}
\caption{\label{alg:cg}Conjugate gradient (CG)}
\begin{algorithmic}[1]
\Require {$n \times n$ symmetric positive definite matrix $A$, length-$n$ vector $b$, and initial approximation $x_1$ to $Ax=b$} 
\Ensure {Approximate solution $x_{i+1}$ to $Ax=b$ with updated residual $r_{i+1}$}
\State {$r_{1}=b-Ax_{1},\, p_{1}=r_{1}$}
\For {$i=1,2,\ldots,$ until convergence}

\State $\alpha_{i}={r}^T_{i}r_{i}/{p}^T_{i}Ap_{i}$
\State $q_{i} = \alpha_{i}p_{i}$

\State $x_{i+1}=x_{i}+q_{i}$
\State $r_{i+1}=r_{i}-Aq_{i}$

\State $\beta_{i}={r}^T_{i+1}r_{i+1}/{r}^T_{i}r_{i}$

\State $p_{i+1}=r_{i+1}+\beta_{i}p_{i}$

\EndFor

\end{algorithmic}
\end{algorithm}

The $s$-step CG method consists of an outer loop, indexed by $k$, which iterates over the blocks of iterations, and an inner loop, which iterates over iterations $j\in\{1,\ldots,s\}$ within a block. For clarity, we globally index iterations by $i\equiv sk+j$. It follows from the properties of CG that for $j\in\{0,\ldots,s\}$, 
\begin{align}
p_{sk+j+1},r_{sk+j+1} &\in {\mathcal{K}}_{s+1}(A,p_{sk+1}) + {\mathcal{K}}_{s}(A,r_{sk+1}) \quad \text{and} \nonumber \\
x_{sk+j+1}-x_{sk+1} &\in {\mathcal{K}}_{s}(A,p_{sk+1}) + {\mathcal{K}}_{s-1}(A,r_{sk+1}). \label{eq:reqsubspaces}
\end{align}
Then the CG vectors for the next $s$ iterations lie in the sum of the column spaces of the matrices 
\begin{align}
\mathcal{P}_{k,s} &= [\rho_0(A)p_{sk+1}, \dots, \rho_s(A)p_{sk+1} ], \text{ with }   \text{span}(\mathcal{P}_{k,s})={\mathcal{K}}_{s+1}(A, p_{sk+1}) \quad\text{and} \nonumber \\
\mathcal{R}_{k,s} &= [\rho_0(A)r_{sk+1}, \dots, \rho_{s-1}(A)r_{sk+1} ],\text{ with }   \text{span}(\mathcal{R}_{k,s})={\mathcal{K}}_{s}(A, r_{sk+1}), \label{eq:cg-krylovbasis}
\end{align}
where $\rho_j(z)$ is a polynomial of degree $j$ satisfying the three-term recurrence
\begin{align}
\rho_0(z) &= 1, \qquad \rho_1(z) = (z-{\theta}_0)\rho_0(z)/{\gamma}_0, \quad \text{and}\nonumber \\
\rho_j(z) &= \left( (z-{\theta}_{j-1})\rho_{j-1}(z)-{\sigma}_{j-2}\rho_{j-2}(z)\right) /{\gamma}_{j-1}.\label{eq:polybasiscg}
\end{align}

We define the $s$-step basis matrix $\mathcal{Y}_{k,s} = [\mathcal{P}_{k,s},\mathcal{R}_{k,s}]$ and define $\underline{\mathcal{Y}}_{k,s}$ to be the same as $\mathcal{Y}_{k,s}$ except with all zeros in columns $s+1$ and $2s+1$. 
Note that the $s$'s in the subscripts are unnecessary here, but will be useful in describing the variable $s$-step CG method later on. 
Under certain assumptions on the nonzero structure of $A$, $\mathcal{Y}_{k,s}$ can be computed in parallel in each outer loop for the same asymptotic latency cost as a single matrix-vector product using the so-called ``matrix powers kernel'' (see~\cite{demmel2007avoiding} for details). Since the columns of $\mathcal{Y}_{k,s}$ satisfy~\eqref{eq:polybasiscg}, we can write 
\begin{equation}
A \underline{\mathcal{Y}}_{k,s} = \mathcal{Y}_{k,s}\mathcal{B}_{k,s},\label{eq:AVVB} 
\end{equation}
where $\mathcal{B}_{k,s}$ is a $(2s+1)\times(2s+1)$ block tridiagonal matrix defined by the 3-term recurrence coefficients in~\eqref{eq:polybasiscg}.

By~\eqref{eq:reqsubspaces}, there exist vectors $p'_{k,j+1}$, $r'_{k,j+1}$, and $x'_{k,j+1}$ that represent the coordinates of the CG iterates $p_{sk+j+1}$, $r_{sk+j+1}$, and $x_{sk+j+1}-x_{sk+1}$, respectively, in $\mathcal{Y}_{k,s}$ for $j\in\{0,\ldots,s\}$. That is,
\begin{align}
p_{sk+j+1}&=\mathcal{Y}_{k,s}p'_{k,j+1},\nonumber\\
r_{sk+j+1}&=\mathcal{Y}_{k,s}r'_{k,j+1}, \quad\text{and} \label{eq:coords}\\
x_{sk+j+1}-x_{sk+1}&=\mathcal{Y}_{k,s}x'_{k,j+1}. \nonumber
\end{align} 
%
%
Therefore, in the inner loop of $s$-step CG, rather than update the CG vectors explicitly, we instead update their coordinates in $\mathcal{Y}_{k,s}$, i.e., for $j\in\{1,\ldots,s\}$,
\begin{align*}
x'_{k,j+1} = x'_{k,j} + \alpha_{sk+j}&p'_{k,j},\quad
r'_{k,j+1} = r'_{k,j} - \alpha_{sk+j}\mathcal{B}_{k,s}p'_{k,j},\\
\text{and} \quad p'_{k,j+1} &= r'_{k,j+1} + \beta_{sk+j}p'_{k,j}.
\end{align*}
Thus the matrix-vector product with $A$ in each iteration becomes a matrix-vector product with the much smaller matrix $\B_{k,s}$. This along with the length-$(2s+1)$ vector updates can be performed locally by each processor in each inner loop iteration. 

We can also reduce communication in computing the inner products. Letting $G_{k,s}=\mathcal{Y}_{k,s}^T\mathcal{Y}_{k,s}$, and using~\eqref{eq:coords} and~\eqref{eq:AVVB}, $\alpha_{sk+j}$ and $\beta_{sk+j}$ can be computed by
\begin{align*}
\alpha_{sk+j} &= \big({r}'^{T}_{k,j}G_{k,s}r'_{k,j}\big) / \big({p}'^{T}_{k,j}G_{k,s}\mathcal{B}_{k,s}p'_{k,j}\big) \quad \text{and}\\
\beta_{sk+j} &= \big({r}'^{T}_{k,j+1}G_{k,s}r'_{k,j+1}\big) / \big({r}'^{T}_{k,j}G_{k,s}{r}'_{k,j}\big).
\end{align*}
The matrix $G_{k,s}$ can be computed with one global synchronization per outer loop iteration, which, in terms of asymptotic parallel latency, costs the same as a single inner product. As $\mathcal{B}_{k,s}$ and $G_{k,s}$ are dimension $(2s+1)\times(2s+1)$, $\alpha_{sk+j}$ and $\beta_{sk+j}$ can be computed locally by each processor within the inner loop.

\begin{algorithm}
\caption{$s$-step conjugate gradient}
\label{alg:scg}
\begin{algorithmic}[1]
\Require {$n \times n$ symmetric positive definite matrix $A$, length-$n$ vector $b$, and initial approximation $x_1$ to $Ax=b$} 
\Ensure {Approximate solution $x_{sk+s+1}$ to $Ax=b$ with updated residual $r_{sk+s+1}$}
\State {$r_{1}=b-Ax_{1},\, p_{1}=r_{1}$}

\For {$k=0,1,\dots,$ until convergence}

\State {Compute $\mathcal{Y}_{k,s}=[\mathcal{P}_{k,s},\, \mathcal{R}_{k,s}]$ according to~\eqref{eq:cg-krylovbasis}}\label{alg:cacg:akx}

\State {Compute $G_{k,s}=\mathcal{Y}_{k,s}^T\mathcal{Y}_{k,s}$} \label{alg:cacg:G}

\State Assemble $\mathcal{B}_{k,s}$ such that~\eqref{eq:AVVB} holds

\State $p'_{k,1}=\left[1,\,0_{1,2s}\right]^{T}$, $r'_{k,1}=\left[0_{1,s+1},\,1,\,0_{1,s-1}\right]^{T}$, $x'_{k,1}=\left[0_{1,2s+1}\right]^T$

\For{$j = 1$ to $s$}

\State $\alpha_{sk+j}=\big({r}'^T_{k,j}G_{k,s}r'_{k,j}\big)/\big({p}'^T_{k,j}G_{k,s}{\mathcal{B}_{k,s}}p'_{k,j}\big)$ \label{alg:$s$-step CG-Method:alpha}
\State $q'_{k,j} = \alpha_{sk+j}p'_{k,j}$

\State $x'_{k,j+1}=x'_{k,j}+q'_{k,j}$ \label{alg:$s$-step CG-Method:x}
\State $r'_{k,j+1}=r'_{k,j}-{\mathcal{B}_{k,s}}q'_{k,j})$ \label{alg:$s$-step CG-Method:r}

\State $\beta_{sk+j}=\big({r}'^T_{k,j+1}G_{k,s}r'_{k,j+1}\big)/\big({r}'^T_{k,j}G_{k,s}r'_{k,j}\big)$ \label{alg:$s$-step CG-Method:beta}

\State $p'_{k,j+1}=r'_{k,j+1}+\beta_{sk+j}p'_{k,j}$ \label{alg:$s$-step CG-Method:p}

\EndFor

\State {Recover iterates $\{p_{sk+s+1},r_{sk+s+1},x_{sk+s+1}\}$ according to \eqref{eq:coords}} \label{alg:cacg:recover}

\EndFor

\end{algorithmic}
\end{algorithm}

%
%
%

\section{Related work}
\label{sec:relatedwork}

In this section, we discussed related work in the areas of variable $s$-step Krylov subspace methods, 
the analysis of maximum attainable accuracy in finite precision CG, and inexact Krylov subspace methods. 

\subsection{Variable $s$-step Krylov subspace methods}
Williams et al.~\cite{williams2014ipdps} use a variable $s$-step BICGSTAB as the coarse 
grid solver in a multigrid method. The technique used here, termed 
``telescoping'', is motivated by the observation than some coarse grid solves 
converge after only a few iterations, whereas other solves take significantly 
longer. By starting with a small $s$ value and gradually increasing it, they 
ensure that $s$-step BICGSTAB with larger $s$ is only used when the solve 
takes enough iterations to amortize the additional costs associated with $s
$-step CG. 

Imberti and Erhel~\cite{imberti2016vary} have recently developed an $s$-step 
GMRES method that allows variable $s$ sizes. They recommend choosing $s_k$ 
values according to the Fibonacci sequence up to some maximum value $s_
{\text{max}}$, i.e., starting with $s_0=s_1=1$ and for $k\geq 2$, $s_k= \min
{(s_{\text{max}}, s_{k-1}+s_{k-2})}$. In this way, the sequence $\{s_i\}$ 
used by Imberti and Erhel is predetermined. 

In contrast, in our approach, the sequence $\{s_i\}$ is dynamically chosen 
based on the 2-norm of the updated residual (which is not necessarily 
monotonically decreasing in CG). In addition, our method is designed such that  
a user-specified accuracy $\veps^*$ can be attained when $\veps^*\geq \veps_{\text{CG}}$, where 
$\veps_{\text{CG}}$ is the accuracy attained by classical CG for the same 
problem. 

\subsection{Maximum attainable accuracy in Krylov subspace methods}
\label{sec:maa}

In both $s$-step and classical variants of Krylov subspace methods, finite precision roundoff error in updates to the approximate solution $x_i$ and the residual $r_i$ in each iteration can cause the updated residual $r_i$ and the true residual $b-Ax_i$ to grow further and further apart as the iterations proceed. If this deviation grows large, it can limit the \emph{maximum attainable accuracy}, i.e., the accuracy with which we can solve $Ax=b$ on a computer with unit round-off $\veps$. 
Analyses of maximum attainable accuracy in CG and other classical KSMs are given by Greenbaum~\cite{greenbaum1995estimating}, van der Vorst and Ye~\cite{van1999residual}, Sleijpen, van der Vorst, and Fokkema~\cite{sleijpen1994bicgstab}, Sleijpen, van der Vorst, and Modersitzki~\cite{sleijpen2001differences}, Bj{\"o}rck, Elfving, and Strako\v{s}~\cite{bjorck1998stability}, and Gutknecht and Strako\v{s}~\cite{gutknecht2000accuracy}.
One important result of these analyses is the insight that loss of accuracy can be caused at a very early stage of the computation, which can not be corrected in later iterations. Analyses of the maximum attainable accuracy in $s$-step CG and the $s$-step biconjugate gradient method (BICG) can be found in~\cite{carson2014residual,carson2015communication}.

\subsection{Inexact Krylov subspace methods}
\label{sec:inexact}

The term ``inexact Krylov subspace methods'' refers to Krylov subspace 
methods in which the matrix-vector products $Av$ are computed inexactly as 
$(A+E)v$ where $E$ represents some error matrix (due to either finite 
precision computation or intentional approximation). It is assumed that all 
other computations (orthogonalization and vector updates) are performed 
exactly. Using an analysis of the deviation of the true and updated 
residuals, it is shown that under these assumptions the size of the 
perturbation term $E$ can be allowed to grow inversely proportional to the 
norm of the updated residual without adversely affecting the attainable 
accuracy. This theory has potential application in mixed-precision Krylov 
subspace methods, as well as in situations where the operator $A$ is 
expensive to apply and can be approximated in order to gain performance.

Much work has been dedicated to inexact Krylov subspace methods. 
In~\cite{simoncinitheory2003}, Simoncini and Szyld provide a general theory 
for inexact variants applicable to both symmetric and nonsymmetric linear 
systems and eigenvalue problems and use this to derive criteria for bounding 
the size of $E$ in such a way that does not diminish the attainable accuracy 
of the method. A similar analysis was given by van de Eschof and Sleijpen~
\cite{vaninexact2004}.
These analyses confirm and improve upon the earlier empirical results for the 
inexact GMRES method of Bouras and Fraysse{\'e}~\cite
{bouras2000relaxationGMRES} and the inexact conjugate gradient method of 
Bouras, Frayss{\'e}, and Giraud~\cite{bouras2000relaxationCG}.

Our work is very closely related to the theory of inexact Krylov subspace 
methods in that our analysis results in a somewhat analogous ``relaxation 
strategy'', where instead of relaxing the accuracy of the matrix-vector 
products we are instead relaxing a constraint on the condition numbers of the 
computed $s$-step basis matrices. In addition, our analysis assumes that all 
parts of the computation are performed in a fixed precision $\veps$.

\section{Variable $s$-step CG}
\label{sec:varscg}

\begin{algorithm}
\caption{Variable $s$-step conjugate gradient}
\label{alg:vscg}
\begin{algorithmic}[1]
\Require {$n \times n$ symmetric positive definite matrix $A$, length-$n$ vector $b$, initial approximation $x_1$ to $Ax=b$, sequence of $s$ values $(s_0,s_1,\ldots)$} 
\Ensure {Approximate solution $x_{m+s_k+1}$ to $Ax=b$ with updated residual $r_{m+s_k+1}$}
\State {$r_{1}=b-Ax_{1},\, p_{1}=r_{1}$}
\State{$m = 0$}

\For {$k=0,1,\dots,$ until convergence}

\State {Compute $\mathcal{Y}_{k,s_k}=[\mathcal{P}_{k,s_k},\, \mathcal{R}_{k,s_k}]$ according to~\eqref{eq:cg-krylovbasis}}

\State {Compute $G_{k,s_k}=\mathcal{Y}_{k,s_k}^T\mathcal{Y}_{k,s_k}$} 

\State Assemble $\mathcal{B}_{k,s_k}$ such that~\eqref{eq:AVVB} holds

\State $p'_{k,1}=\left[1,\,0_{1,2s_k}\right]^{T}$, $r'_{k,1}=\left[0_{1,s_k+1},\,1,\,0_{1,s_k-1}\right]^{T}$, $x'_{k,1}=\left[0_{1,2s_k+1}\right]^T$

\For{$j = 1$ to $s_k$}

\State $\alpha_{m+j}=\big({r}'^T_{k,j}G_{k,s_k}r'_{k,j}\big)/\big({p}'^T_{k,j}G_{k,s_k}{\mathcal{B}_{k,s_k}}p'_{k,j}\big)$ 
\State $q'_{k,j} = \alpha_{m+j}p'_{k,j}$

\State $x'_{k,j+1}=x'_{k,j}+q'_{k,j}$
\State $r'_{k,j+1}=r'_{k,j}-{\mathcal{B}_{k,s_k}}q'_{k,j})$ 

\State $\beta_{m+j}=\big({r}'^T_{k,j+1}G_{k,s_k}r'_{k,j+1}\big)/\big({r}'^T_{k,j}G_{k,s_k}r'_{k,j}\big)$ 

\State $p'_{k,j+1}=r'_{k,j+1}+\beta_{m+j}p'_{k,j}$ 

\EndFor

\State {Recover iterates $\{p_{m+s_k+1},r_{m+s_k+1},x_{m+s_k+1}\}$ according to \eqref{eq:coords}}  
\State {$m=m+s_k$}

\EndFor

\end{algorithmic}
\end{algorithm}

As discussed, the con\-di\-tion\-ing of the $s$-step basis matrices in $s$-step Kry\-lov sub\-space meth\-ods affects the 
rate of convergence and attainable accuracy. 
In the $s$-step CG method (Algorithm~\ref{alg:scg}), the Krylov subspace basis is computed in blocks of size $s$ in each outer loop. 
This method can be generalized to allow the blocks to be of varying size; in other words, a different $s$ value can be used in each outer 
loop iteration. We denote the block size in outer iteration $k$ by $s_k$, where $1\leq s_k \leq s_{\text{max}}$ for some maximum $s$ value 
$s_{\text{max}}$. We call this a \emph{variable $s$-step Krylov subspace method}. A general variable $s$-step CG method is shown in Algorithm~\ref{alg:vscg}. 

In this section, we derive a variable $s$-step CG method which automatically determines $s_k$ in each outer loop such that a user-specified accuracy can be attained, which we call the adaptive $s$-step CG method. 
Our analysis is based on the maximum attainable accuracy analysis for $s$-step CG in~\cite{carson2014residual}, which shows that the attainable accuracy of $s$-step CG depends on the condition numbers of the $s$-step basis matrices $\y_{k,s}$ computed in each outer loop. In summary, we show that if in outer loop $k$, $s_k$ is selected such that the condition number of the basis matrix $\y_{k,s_k}$ is inversely proportional to the maximum 2-norm of the residual vectors computed within outer loop $k$, then the approximate solution produced by the variable $s$-step CG method can be as accurate as the approximate solution produced by the classical CG method. Further, our method allows for the user to specify the desired accuracy, so in the case that a less accurate solution is required, our method automatically adjusts to allow for higher $s_k$ values. This effectively exploits the tradeoff between accuracy and performance in $s$-step Krylov subspace methods.   

In this section we derive our adaptive $s$-step approach, which stems from a bound on the 
growth in outer loop $k$ of the gap between the true and updated residuals in finite precision. 

\subsection{A constraint on basis matrix condition number}
\label{sec:theory}

In the remainder of the paper, hats denote quantities computed in finite precision, $\veps$ denotes the unit round-off, $\norm{\cdot}$ denotes the 2-norm, and $\kappa(A)$ denotes the 2-norm condition number $\norm{A^{-1}}\norm{A}$. 
To simplify the analysis, we drop all $O(\veps^2)$ terms (and higher order terms in $\veps$). 
Writing the true residual $b-A\xhat_i=b-A\xhat_i-\rhat_i + \rhat_i$, we can write the bound
\[
\norm{b-A\xhat_i} \leq \norm{b-A\xhat_i-\rhat_i}+\norm{\rhat_i}.
\]
Then as $\norm{\rhat_i}\rightarrow 0$ (i.e., as the method converges), we have $\norm{b-A\xhat_i}\approx \norm{b-A\xhat_i-\rhat_i}$.  
The size of the \emph{residual gap}, i.e., $\norm{b-A\xhat_i -\rhat_i}$, therefore determines the attainable accuracy of the method (see the related references in section~\ref{sec:maa}). 

We begin by considering the rounding errors made in the variable $s$-step CG method (Algorithm~\ref{alg:vscg}). 
In outer loop $k$ and inner loops $j\in\{1,\ldots,s_k\}$ of finite precision variable $s$-step CG, using standard models of floating point error (e.g.,~\cite[\S2.4]{golub1996matrix}) we have
\begin{align}
A\uyh_{k,s_k} &= \yh_{k,s_k} \B_{k,s_k} - E_k, \label{E}\\ 
&\text{where}\quad \norm{E_k} \leq \veps \left( (3+N_A)\norm{A}\norm{\uyh_{k,s_k}} + 4\norm{\yh_{k,s_k}}\norm{\B_{k,s_k}} \right), \nonumber\\
\xhat'_{k,j+1} &= \xhat'_{k,j} + \qhat'_{k,j} + \xi_{k,j+1},\label{xp} \\ 
&\text{where}\quad \norm{\xi_{k,j+1}}\leq \veps\norm{\xhat'_{k,j+1}},\nonumber\\
\rhat'_{k,j+1} &= \rhat'_{k,j}-\B_{k,s_k} \qhat'_{k,j} + \eta_{k,j+1}, \label{rp}\\ 
&\text{where}\quad \norm{\eta_{k,j+1}}\leq \veps (\norm{\rhat'_{k,j+1}}+3\norm{\B_{k,s_k}}\norm{\qhat'_{k,j}}),\nonumber\\
\xhat_{m+j+1} &= \yh_{k,s_k}\xhat'_{k,j+1} + \xhat_{m+1} + \phi_{m+j+1},\label{x}\\ 
&\text{where}\quad \norm{\phi_{m+j+1}}\leq \veps (\norm{\xhat_{m+j+1}} + 3j\norm{\yh_{k,s_k}}\norm{\xhat'_{k,j+1}}),\text{ and}\nonumber\\
\rhat_{m+j+1}&= \yh_{k,s_k} \rhat'_{k,j+1} + \psi_{m+j+1},\label{r}, \\ 
&\text{where}\quad \norm{\psi_{m+j+1}}\leq 3j\veps\norm{\yh_{k,s_k}}\norm{\rhat'_{k,j+1}}\nonumber, 
\end{align}
where $N_A$ denotes the maximum number of nonzeros per row in $A$ and $m=\sum_{i=0}^{k-1} s_i$. 
Let $\delta_i \equiv b-A\xhat_i - \rhat_i$ denote the residual gap in iteration $i$. Similar to the analysis in~\cite{carson2014residual}, for variable $s$-step CG we can write the growth of the residual gap in iteration $m+j+1$ in outer loop $k$ as 
\begin{equation*}
\delta_{m+j+1} - \delta_{m+1} = - \sum_{i=1}^j \left(A\uyh_{k,s_k} \xi_{k,i+1}+\yh_{k,s_k} \eta_{k,i+1}  -E_k\qhat'_{k,i} \right) - A\phi_{m+j+1} - \psi_{m+j+1}.
\end{equation*}
%
%

Our goal is now to manipulate a bound on this quantity in order to derive a method for determining the largest $s_k$ value we can use in outer loop $k$ such that the desired accuracy is attained. Using standard techniques, we have the norm-wise bound
\begin{align*}
\norm{\delta_{m+j+1}-\delta_{m+1}} 
&\leq \veps\norm{A}\norm{\uyh_{k,s_k}}\sum_{i=1}^{j} \norm{\xhat'_{k,i+1}} + \veps \norm{\yh_{k,s_k}}\sum_{i=1}^{j} \norm{\rhat'_{k,i+1}} \\
&\phantom{\leq} + \veps\left( (3+N_A)\norm{A}\norm{\uyh_{k,s_k}} + 7 \norm{\yh_{k,s_k}}\norm{\B_{k,s_k}} \right) \sum_{i=1}^{j} \norm{\qhat'_{k,i}} \\
&\phantom{\leq} + \veps\norm{A}\norm{\xhat_{m+j+1}} + 3j\veps \norm{A}\norm{\yh_{k,s_k}} \norm{\xhat'_{k,j+1}} + 3j\veps \norm{\yh_{k,s_k}} \norm{\rhat'_{k,j+1}}\\
&\leq \veps \left( (7+2N_A)\norm{A}\norm{\uyh_{k,s_k}} + 14 \norm{\yh_{k,s_k}}\norm{\B_{k,s_k}} \right)\sum_{i=1}^{j} \norm{\xhat'_{k,i+1}} \\
&\phantom{\leq} + \veps\norm{\yh_{k,s_k}}\sum_{i=1}^{j} \norm{\rhat'_{k,i+1}}\\
&\phantom{\leq} + \hspace{-1pt}\veps\norm{A}\norm{\xhat_{m+j+1}} \hspace{-1pt}+\hspace{-1pt} 3j\veps \norm{A}\norm{\yh_{k,s_k}} \norm{\xhat'_{k,j+1}} \hspace{-1pt}+\hspace{-1pt} 3j\veps \norm{\yh_{k,s_k}} \norm{\rhat'_{k,j+1}},
\end{align*}
where we have used $\sum_{i=1}^{j}\norm{\qhat'_{k,i}}\leq (2+\veps)\sum_{i=1}^{j}\norm{\xhat'_{k,i+1}}$.
To simplify the notation, we let $\tau_k = \norm{\B_{k,s_k}}/\norm{A}$ and $N_{k}=7+2N_A+14\tau_k$. We note that $\norm{\B_{k,s_k}}$ depends on the polynomial basis used and is not expected to be too large; e.g., for the monomial basis, $\norm{\B_{k,s_k}}=1$.
The above bound can then be written
\begin{align}
\norm{\delta_{m+j+1}-\delta_{m+1}} &\leq \veps N_{k}\norm{A}\norm{\yh_{k,s_k}} \sum_{i=1}^{j} \norm{\xhat'_{k,i+1}} + \veps\norm{\yh_{k,s_k}} \sum_{i=1}^{j} \norm{\rhat'_{k,i+1}} \label{b1}\\
&\phantom{\leq} +\hspace{-1pt} \veps\norm{A}\norm{\xhat_{m+j+1}} \hspace{-1pt}+\hspace{-1pt} 3j\veps \norm{A}\norm{\yh_{k,s_k}} \norm{\xhat'_{k,j+1}} \hspace{-1pt}+\hspace{-1pt} 3j\veps \norm{\yh_{k,s_k}} \norm{\rhat'_{k,j+1}}.\nonumber
\end{align}

Note that assuming $\yh_{k,s_k}$ is full rank, we have
\begin{align}
\norm{\yh_{k,s_k}}\norm{\rhat'_{k,i+1}} &= \norm{\yh_{k,s_k}}\norm{\yh_{k,s_k}^{+} \yh_{k,s_k} \rhat'_{k,i+1}}\nonumber\\
&\leq \kappa(\yh_{k,s_k}) \norm{\yh_{k,s_k} \rhat_{m+i+1}-\psi_{m+i+1}}\nonumber\\
&=\kappa(\yh_{k,s_k})\norm{\rhat_{m+i+1}}+O(\veps), \quad\text{and}, \label{yr} \\
\norm{\yh_{k,s_k}}\norm{\xhat'_{k,i+1}} &= \norm{\yh_{k,s_k}}\norm{\yh_{k,s_k}^{+} \yh_{k,s_k} \xhat'_{k,i+1}} \nonumber \\
&\leq \kappa(\yh_{k,s_k}) \norm{\yh_{k,s_k} \xhat'_{k,i+1}}. \label{yx}\\
\end{align}
To bound $\norm{\yh_{k,s_k} \xhat'_{k,i+1}}$ in terms of the size of the residuals, notice that by~\eqref{x}, 
\begin{align*}
\yh_{k,s_k} \xhat'_{k,i+1} &= \xhat_{m+i+1} - \xhat_{m+1} - \phi_{m+i+1}\\
&= (x-\xhat_{m+1}) - (x-\xhat_{m+i+1})-\phi_{m+i+1}\\
&= A^{-1}(b-A\xhat_{m+1}) - A^{-1} (b-A\xhat_{m+i+1}) - \phi_{m+i+1}\\
&= A^{-1}(\rhat_{m+1}+\delta_{m+1}) - A^{-1}(\rhat_{m+i+1}+\delta_{m+i+1})-\phi_{m+i+1}\\
&= A^{-1}(\rhat_{m+1}-\rhat_{m+i+1}) + A^{-1}(\delta_{m+1}-\delta_{m+i+1}) - \phi_{m+i+1},
\end{align*}
so
\begin{equation*}
\norm{\yh_{k,s_k}\xhat'_{k,i+1}}\leq \norm{A^{-1}}(\norm{\rhat_{m+1}}+ \norm{\rhat_{m+i+1}}) + O(\veps). 
\end{equation*}
Thus together with~\eqref{yx} we can write
\begin{equation*}
\norm{\yh_{k,s_k}} \norm{\xhat'_{k,i+1}} \leq \kappa(\yh_{k,s_k})\norm{A^{-1}}(\norm{\rhat_{m+1}}+\norm{\rhat_{m+i+1}})+O(\veps).
\end{equation*}

This allows us to write the bound~\eqref{b1} in the form
\begin{align*}
\norm{\delta_{m+j+1}-\delta_{m+1}} 
&\leq 2\veps N_k\kappa(A)\kappa(\yh_{k,s_k})\sum_{i=1}^{j} \norm{\rhat_{m+i+1}}\\
&\phantom{\leq}+\veps\left(N_k j\kappa(A)\kappa(\yh_{k,s_k})+3j\kappa(A)\kappa(\yh_{k,s_k}) \right) \norm{\rhat_{m+1}}\\
&\phantom{\leq}+\veps\left( 3j\kappa(A)\kappa(\yh_{k,s_k})+3j\kappa(\yh_{k,s_k}) \right) \norm{\rhat_{m+j+1}}\\
&\phantom{\leq}+\veps\norm{A}\norm{\xhat_{m+j+1}}.
\end{align*}
Using $x-\xhat_{m+j+1} = A^{-1} \rhat_{m+j+1} + \delta_{m+j+1}$, we have
%
%
\[
\norm{\xhat_{m+j+1}-x}\leq \norm{A^{-1}}\norm{\rhat_{m+j+1}}+O(\veps),
\]
and writing $\xhat_{m+j+1}=\xhat_{m+j+1}-x+x$, we have
\begin{align*}
\norm{\delta_{m+j+1}-\delta_{m+1}} 
%
%
&\leq 2N_k \veps\kappa(A)\kappa(\yh_{k,s_k})\left( j (\norm{\rhat_{m+1}}+\norm{\rhat_{m+j+1}}) + \sum_{i=1}^{j}\norm{\rhat_{m+i+1}} \right)\\
&\phantom{\leq}+\veps\norm{A}\norm{x}.
\end{align*}
%
Then
\begin{equation*}
\norm{\delta_{m+j+1}-\delta_{m+1}} \leq 6N_kj\veps\kappa(A)\kappa(\yh_{k,s_k}) \left( \max_{\ell\in\{0,\ldots,j\}} \norm{\rhat_{m+\ell+1}}\right) +\veps\norm{A}\norm{x},
\end{equation*}
and letting $c_{k,j}=6N_k j\kappa(A)$, this bound can be written
\begin{equation}
\norm{\delta_{m+j+1}-\delta_{m+1}} \leq c_{k,j}\veps \kappa(\yh_{k,s_k}) \left(\max_{\ell\in\{0,\ldots,j\}}\norm{\rhat_{m+\ell+1}} \right) + \veps\norm{A}\norm{x}.
\label{fb}
\end{equation}
This gives a norm-wise bound on the growth of the residual gap due to finite precision errors made in outer iteration $k$.

Letting $m_\ell = \sum_{i=0}^{\ell-1} s_i$, notice that we have
\[
\norm{\delta_{m_k+s_k+1}}\leq \norm{\delta_1} + \sum_{\ell=0}^{k} \norm{\delta_{m_\ell+s_\ell+1}- \delta_{m_\ell+1}}.
\]
Suppose that at the end of the iterations we want the size of the residual gap to be on the order $\veps^*$ (where $\veps^* \gtrsim O(\veps)\norm{A}\norm{x}$, as this is the best we could expect). 
Then assuming the number of outer iterations $k$ is not too high,  
this can be accomplished by requiring that in each outer loop, for 
$j\in \{1,\ldots, s_k\}$,
\[
\norm{\delta_{m+j+1}-\delta_{m+1}} \leq \veps^*.
\]
From~\eqref{fb}, this means that we must have, for all inner iterations $j\in\{1,\ldots, s_k\}$,
\begin{equation}
\kappa(\yh_{k,s_k}) \leq \frac{\veps^*}{c_{k,j}\veps\max_{\ell\in\{0,\ldots,j\}}\norm{\rhat_{m+\ell+1}}}.
\label{yhk}
\end{equation}
The bound~\eqref{yhk} tells us that as the 2-norm of the residual decreases, we can tolerate a more ill-conditioned basis matrix $\yh_{k,s}$ without detriment to the attainable accuracy. Since $\kappa(\yh_{k,s})$ grows with increasing $s_k$, this suggests that $s_k$ can be allowed to increase as the method converges at a rate proportional to the inverse of the residual 2-norm. 
This naturally gives a relaxation strategy that is somewhat analogous to those derived for inexact Krylov subspace methods (see section~\ref{sec:inexact}). We discuss the implementation of this strategy in the following subsection.

\subsection{The adaptive $s$-step CG method}
\label{sec:algorithm}

In this section, we propose a variable $s$-step CG method that makes use of the constraint~\eqref{yhk} in determining $s_k$ in each outer loop $k$, which we call the \emph{adaptive} $s$-step CG method. Our approach uses the largest $s_k$ value possible in each outer iteration $k$ such that~\eqref{yhk} is satisfied up to some $s_{\text{max}}$ set by the user (e.g., $s_{\text{max}}$ could be selected by auto-tuning to find the $s$ value that minimizes the time per iteration in $s$-step CG).  Since a variable $s$-step CG method produces the same iterates as classical CG when $s_k=1$ for all outer loops $k$, we expect that our adaptive $s$-step CG method can attain accuracy $\veps^*$ when $\veps^*\geq \veps_{\text{CG}}$, where $\veps_{\text{CG}}$ is the attainable accuracy of classical CG. When the right-hand-side of~\eqref{yhk} is less than one (due to, e.g., $\veps^*$ being set too small), the inequality~\eqref{yhk} cannot be satisfied. In this case, we default to using $s_k=1$ (i.e., we perform an iteration of classical CG). 

In CG, it is the $A$-norm of the error rather than the $2$-norm of the residual that is minimized, which means that in some inner loop iteration $j$, we can have $\norm{\rhat_{m+j+1}}> \norm{\rhat_{m+1}}$. Since $\norm{\rhat_{m+i+1}}$ for $i>1$ are not known at the start of the outer loop iteration, we must determine $s_k$ online, i.e., during the iterations. We therefore first construct a basis matrix that satisfies~\eqref{yhk} under the assumption that $\norm{\rhat_{m+i+1}} \leq \norm{\rhat_{m+1}}$. Within each inner loop, we can perform an inexpensive check to ensure that~\eqref{yhk} is still satisfied after computing the next residual; if not, we break from the inner loop, discard the current basis matrix, and begin a new outer loop iteration. We elaborate on this approach in the remainder of this section. 
An optimized high-performance implementation of our variable $s$-step 
approach is outside the scope of the present work, although we will briefly 
discuss how one might implement this strategy efficiently such that each outer loop still requires only a single global synchronization and the number of wasted matrix-vector products is kept to a minimum.

The bound~\eqref{yhk} indicates that when $\yh_{k,s_k}$ is constructed at the beginning of outer loop $k$, we should enforce  
\begin{equation}
\kappa(\yh_{k,s_k}) \leq \frac{\veps^*}{c_{k,s_{\text{max}}}\veps\norm{\rhat_{m+1}}}.
\label{yhk1}
\end{equation}
%
Therefore we first construct an $\bar{s}_k$-step basis matrix $\yh_{k,\bar{s}_k}$ using some $\bar{s}_k \leq s_{\text{max}}$, compute the Gram matrix $\Gh_{k,\bar{s}_k} = \yh_{k,\bar{s}_k}^T \yh_{k,\bar{s}_k}$, and then use  
the condition numbers of the leading principle submatrices of the blocks of $\Gh_{k,\bar{s}_k}$ to compute the condition numbers of the leading columns of the blocks of $\yh_{k,\bar{s}_k}$. To be concrete, for $i\in \{1,2,\ldots,\bar{s}_k\}$, we have (using MATLAB notation) 
\begin{equation*}
\yh_{k,i} = \left[\yh_{k,\bar{s}_k}({:},1{:}i{+}1), \yh_{k,\bar{s}_k}({:},i_s{:}i_s{+}i{-}1)\right]
\end{equation*}
and
\begin{equation*}
\Gh_{k,i} = 
\begin{bmatrix}
\Gh_{k,\bar{s}_k}(1{:}i{+}1,1{:}i{+}1) & \Gh_{k,\bar{s}_k}(1{:}i{+}1, i_s{:}i_s{+}i{-}1) \\
\Gh_{k,\bar{s}_k}(i_s{:}i_s{+}i{-}1,1{:}i{+}1) & \Gh_{k,\bar{s}_k}(i_s{:}i_s{+}i{-}1,i_s{:}i_s{+}i{-}1)
\end{bmatrix},
\end{equation*}
where $i_s = \bar{s}_k+2$.
Then since $\sqrt{\kappa(\Gh_{k,i})} = \kappa(\yh_{k,i})$, 
we can easily compute the condition number of the $i$-step basis matrix $\yh_{k,i}$, i.e., the basis matrix needed to perform $i$ inner loop iterations, for $i\in\{1,\ldots,\bar{s}_k\}$. We therefore let $\tilde{s}_k\in \{1,\ldots,\bar{s}_k\}$ be the largest value possible such that 
\begin{equation}
\sqrt{\kappa(\Gh_{k,\tilde{s}_k})} \leq \frac{\veps^*}{c_{k,s_{\text{max}}}\veps\norm{\rhat_{m+1}}}
\label{condy}
\end{equation}
holds. 

%
 
To handle the case where $\norm{\rhat_{m+j+1}}>\norm{\rhat_{m+1}}$ for some inner iteration $j$, we can add a check within the inner loop iterations: after computing $\rhat'_{k,j+1}$, we check if
\begin{equation}
\kappa(\yh_{k,\tilde{s}_k}) \leq \frac{\veps^*}{c_{k,j}\veps\norm{\rhat_{m+j+1}}} \approx \frac{\veps^*}{c_{k,j}\veps\left(\rhat_{k,j+1}^{'T} \Gh_{k,\tilde{s}_k} \rhat'_{k,j+1} \right)^{1/2}}
\label{check}
\end{equation}
is still satisfied. Note that  
since $\norm{\rhat_{m+j+1}} = \left(\rhat_{k,j+1}^{'T} \Gh_{k,\tilde{s}_k} \rhat'_{k,j+1} \right)^{1/2} + O(\veps)$, performing this check 
does not require any additional communication.
If~\eqref{check} is satisfied for all $j\in \{1,\ldots,\tilde{s}_k\}$, we will perform all $s_k=\tilde{s}_k$ inner loop iterations. If this is not satisfied for some $j$, we break from the current inner loop iteration and begin the next outer loop; i.e., we only perform $s_k=j$ inner loop iterations in outer loop $k$. 

Although performing the checks for determining when to break from the current outer loop do not incur any additional communication cost in terms of the number of messages sent between processors, there is, however, potential wasted effort in terms of computation and the number of words moved when $\bar{s}_k > \tilde{s}_k$ and/or $\tilde{s}_k>s_k$, since in this case we have computed more basis vectors than we need to perform $s_k$ inner loop updates.  
One way to mitigate this is to only allow 
$\bar{s}_k$ to grow by at most $f$ vectors in each outer loop, i.e., $\bar{s}_k\leq s_
{k-1}+f \leq s_{\text{max}}$, where $f$ is some small positive integer (e.g., $f=1$ or $f=2$).

The approach outlined above is summarized in Algorithm~\ref{alg:vscg2}. Note that in Algorithm~\ref{alg:vscg2}, we have specified that the function $c_{k,j}$ can be input by the user. Because our bounds are not tight, in many cases using $c_{k,j}=6N_kj\kappa(A)$ (from the bound~\eqref{yhk}) overestimates the error, which results in the use of smaller $s_k$ values than necessary. 

Note that although there is no expectation that $\norm{\rhat_{m+j+1}}$ decreases between iterations in CG, as long as the method 
is converging, we still expect the residual norm to trend downward. Therefore one could make the simplifying assumption that 
$\norm{\rhat_{m+1}}\geq \norm{\rhat_{m+j+1}}$ for $j\in\{1,\ldots,\tilde{s}_k\}$. In this case, the value $s_k=\tilde{s}_k$ can be determined at the beginning of outer loop $k$ before performing any inner iterations, and so the checks~\eqref{check} are unnecessary. This simplification comes at the potential cost 
of reduced reliability of the adaptive $s$-step method, although in our experimenting we have observed that taking $s_k=\tilde{s}_k$ usually works well even when the residual norm is nonmonotonic. 

Because our strategy is based on a loose upper bound on the 
deviation of residuals, we do not require particularly accurate estimates of 
$\kappa(\yh_{k,i})$. Therefore alternative methods for inexpensively approximating this quantity 
could also be used in practice. 

\begin{algorithm}
\caption{Adaptive $s$-step conjugate gradient}
\label{alg:vscg2}
\begin{algorithmic}[1]
\Require {$n \times n$ symmetric positive definite matrix $A$, length-$n$ vector $b$, initial approximation $x_1$ to $Ax=b$, maximum $s$ value $s_{\text{max}}$, initial $s$ value $\bar{s}_0$, maximum basis growth factor $f$, desired convergence tolerance $\veps^*$, function $c_{k,j}$} 
\Ensure {Approximate solution $x_{m+s_k+1}$ to $Ax=b$ with updated residual $r_{m+s_k+1}$}
\State {$r_{1}=b-Ax_{1},\, p_{1}=r_{1}$}
\State{$m = 0$}
\For {$k=0,1,\dots,$ until convergence}
\If {$k\neq 0$}
\State {$\bar{s}_k = \min(s_{k-1}+f, s_{\text{max}})$}
\EndIf
\State {Compute $\bar{s}_k$-step basis matrix $\mathcal{Y}_{k,\bar{s}_k}=[\mathcal{P}_{k,\bar{s}_k},\, \mathcal{R}_{k,\bar{s}_k}]$ according to~\eqref{eq:cg-krylovbasis}}
\State {Compute $G_{k,\bar{s}_k}=\mathcal{Y}_{k,\bar{s}_k}^T\mathcal{Y}_{k,\bar{s}_k}$} 
\State {Determine $\tilde{s}_k$ by~\eqref{condy}; assemble $\mathcal{Y}_{k,\tilde{s}_k}$ and $G_{k,\tilde{s}_k}$}
\State {Store estimate $\gamma \approx \kappa(\yh_{k,\tilde{s}_k})$}
\State Assemble $\mathcal{B}_{k,\tilde{s}_k}$ such that~\eqref{eq:AVVB} holds

\State $p'_{k,1}=\left[1,\,0_{1,2\tilde{s}_k}\right]^{T}$, $r'_{k,1}=\left[0_{1,\tilde{s}_k+1},\,1,\,0_{1,\tilde{s}_k-1}\right]^{T}$, $x'_{k,1}=\left[0_{1,2\tilde{s}_k+1}\right]^T$
\For{$j = 1$ to $\tilde{s}_k$}
\State $s_k=j$
\State $\alpha_{m+j}=\big({r}'^T_{k,j}G_{k,\tilde{s}_k}r'_{k,j}\big)/\big({p}'^T_{k,j}G_{k,\tilde{s}_k}\mathcal{B}_{k,\tilde{s}_k}p'_{k,j}\big)$ 
\State $q'_{k,j} = \alpha_{m+j}p'_{k,j}$

\State $x'_{k,j+1}=x'_{k,j}+q'_{k,j}$
\State $r'_{k,j+1}=r'_{k,j}-{\mathcal{B}_{k,\tilde{s}_k}}q'_{k,j})$ 

\State $\beta_{m+j}=\big({r}'^T_{k,j+1}G_{k,\tilde{s}_k}r'_{k,j+1}\big)/\big({r}'^T_{k,j}G_{k,\tilde{s}_k}r'_{k,j}\big)$ 

\State $p'_{k,j+1}=r'_{k,j+1}+\beta_{m+j}p'_{k,j}$ 

\If {$\gamma \geq \veps^*/ c_{k,j} \veps \left(\rhat_{k,j+1}^{'T} G_{k,\tilde{s}_k} \rhat'_{k,j+1} \right)^{1/2}$}
\State {break from inner loop}
\EndIf

\EndFor

\State {Recover iterates $\{p_{m+s_k+1},r_{m+s_k+1},x_{m+s_k+1}\}$ according to~\eqref{eq:coords}}
\State {$m=m+s_k$}

\EndFor

\end{algorithmic}
\end{algorithm}


\section{Numerical experiments}
\label{sec:experiments}

In this section, we test the adaptive $s$-step CG method (Algorithm~\ref{alg:vscg2}) on a number of small examples from the University of Florida Sparse Matrix Collection~\cite{davis2011university} using MATLAB. The test examples were chosen to exemplify various convergence trajectories of CG, which allows us to demonstrate the validity of our relaxation condition in which $s_k$ depends on the rate of convergence. The test problems used here are small, so we would not use (fixed or variable) $s$-step methods here in practice; nevertheless, they still serve to exhibit the numerical behavior of the methods. 

All tests were run in double precision, i.e., $\veps \approx 2^{-53}$, and 
use a right-hand side $b$ with entries $1/\sqrt{n}$ where $n$ is the dimension of $A$.
We use matrix equilibration, i.e., $A\gets D^{-1/2}AD^{-1/2}$ where $D$ is a diagonal matrix of the largest entries in each row of $A$.  Properties of the test matrices (post-equilibration) are shown in Table~\ref{mats}.

\begin{figure}[ht]
	\centering
			\captionof{figure}{Convergence of the 2-norm of the true residual for the matrix gr\_30\_30. Plots on the left use $\veps^*=$3.4e-14 and 
	plots on the right use $\veps^*=$1e-6. Plots are shown for three $s$ (or $s_{\text{max}}$) values: $4$ (top), $8$ (middle), and $10$ (bottom).  
	Iterations where communication takes place are noted with markers (red circles for $s$-step CG, blue stars for variable $s$-step CG, and black dots for classical CG). The horizontal dashed line shows the requested accuracy $\veps^*$. }
	\vspace{.2cm}
	\includegraphics[width=2.15in]{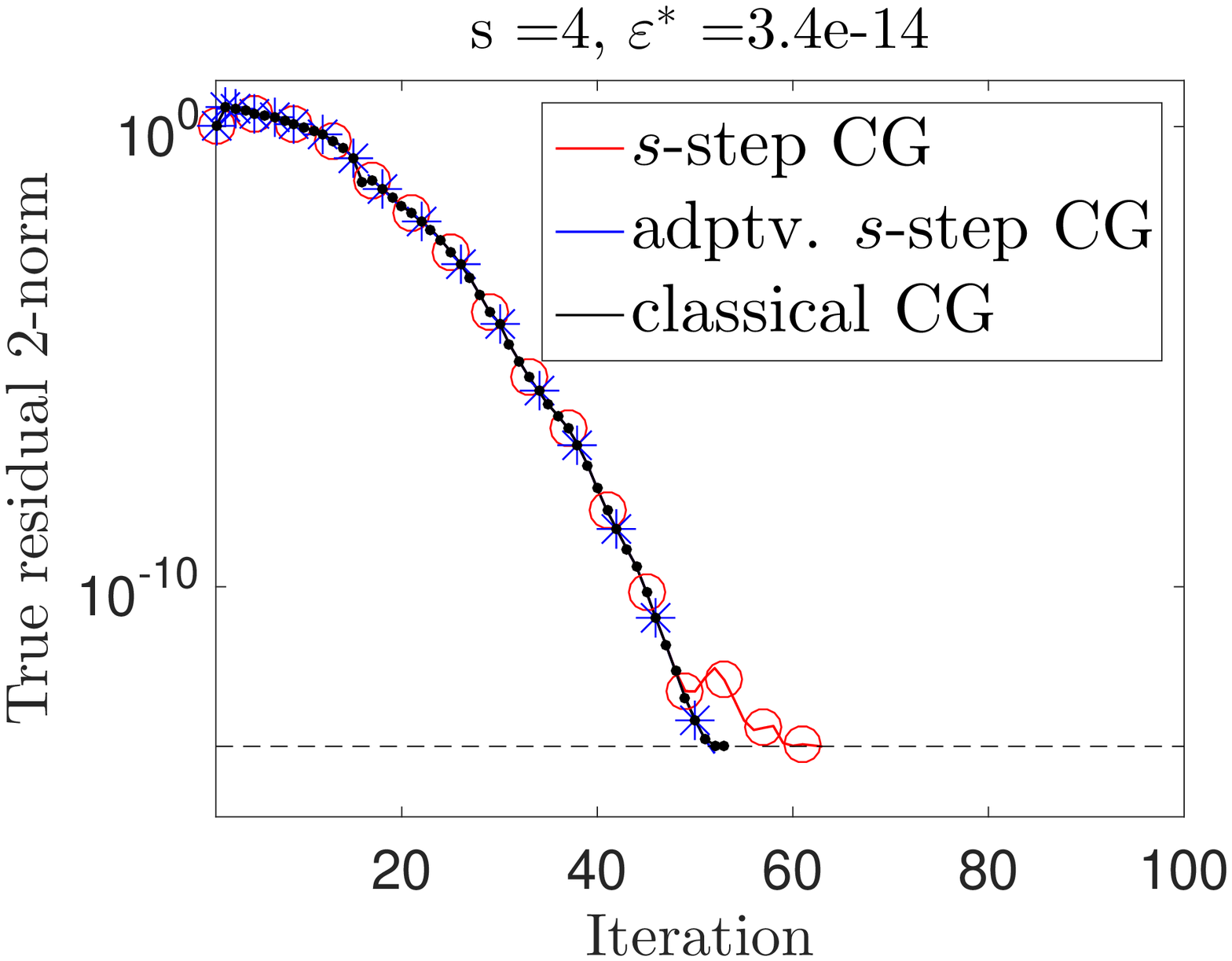}
	\hspace{5mm}
		\includegraphics[width=2.15in]{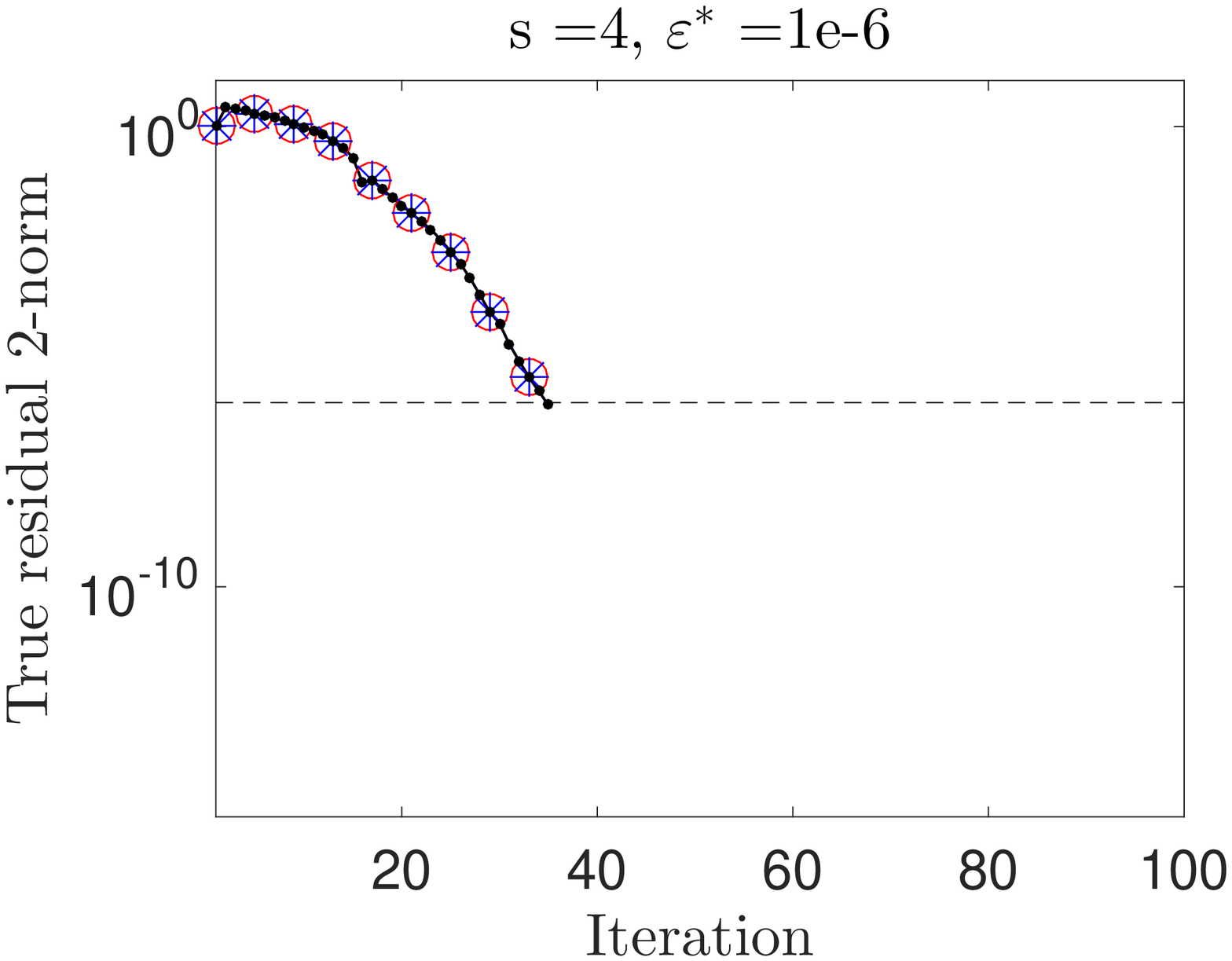}\\
		\vspace{1mm}
	\includegraphics[width=2.15in]{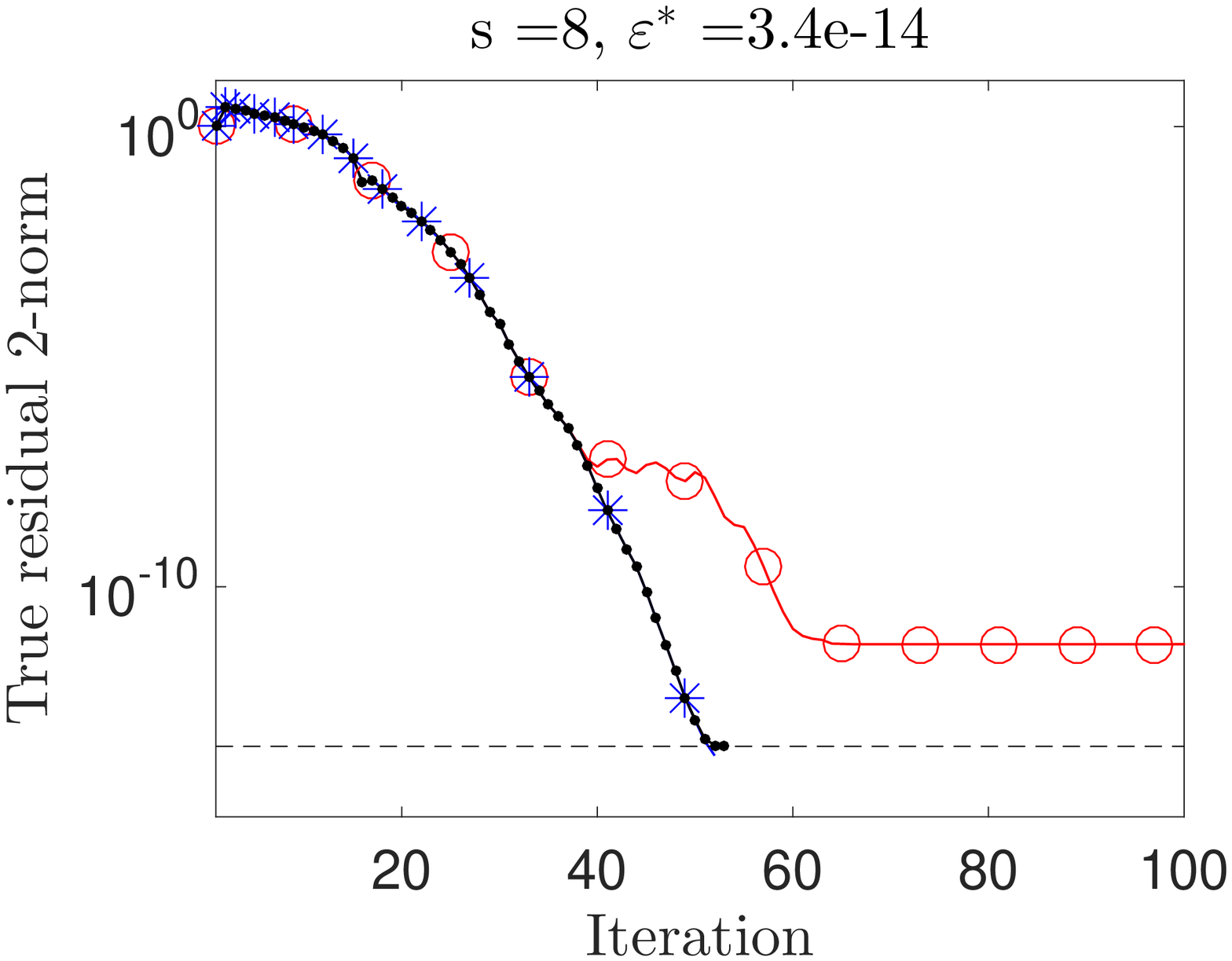}
	\hspace{5mm}
		\includegraphics[width=2.15in]{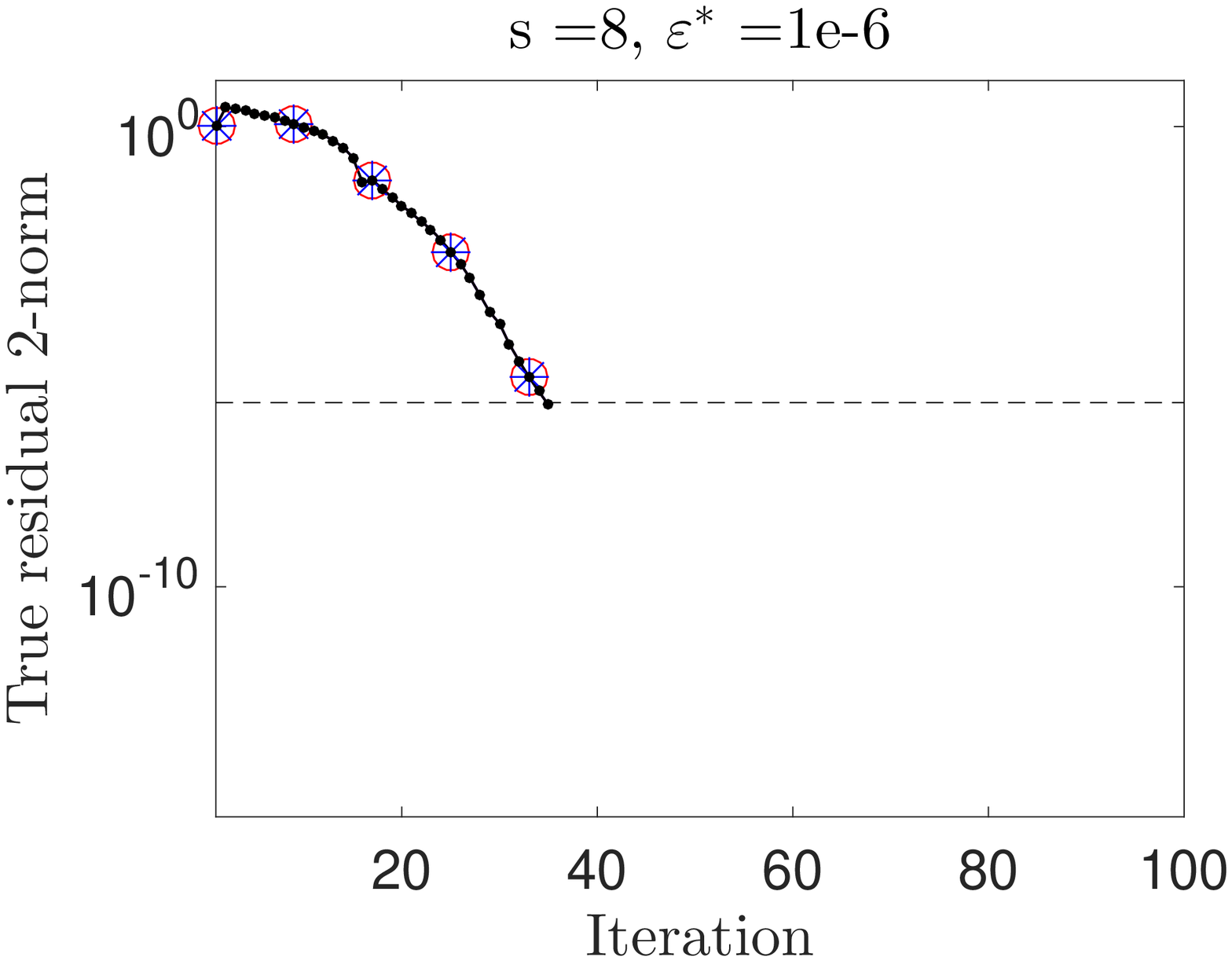}\\
			\vspace{1mm}
			\includegraphics[width=2.15in]{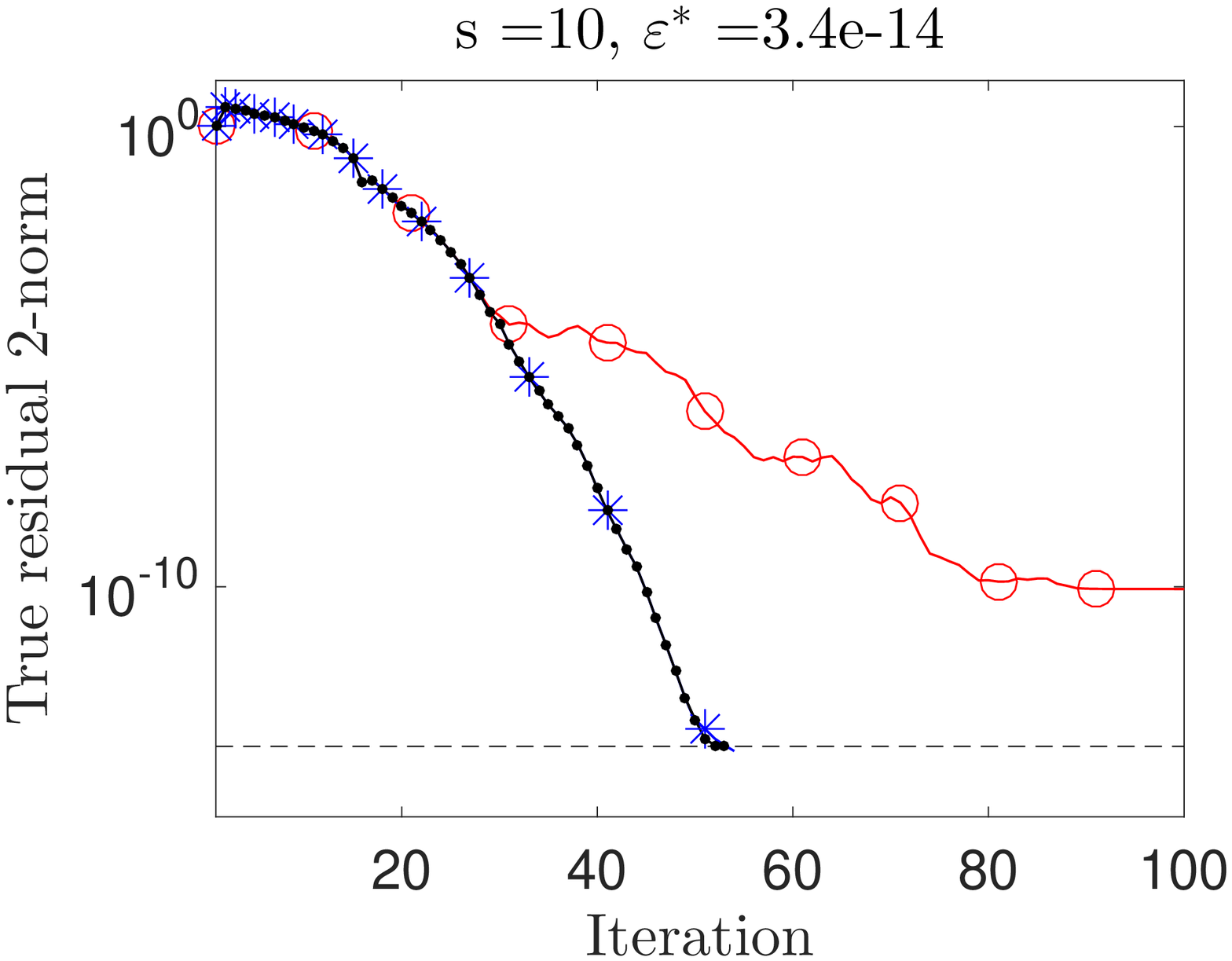}
	\hspace{5mm}
		\includegraphics[width=2.15in]{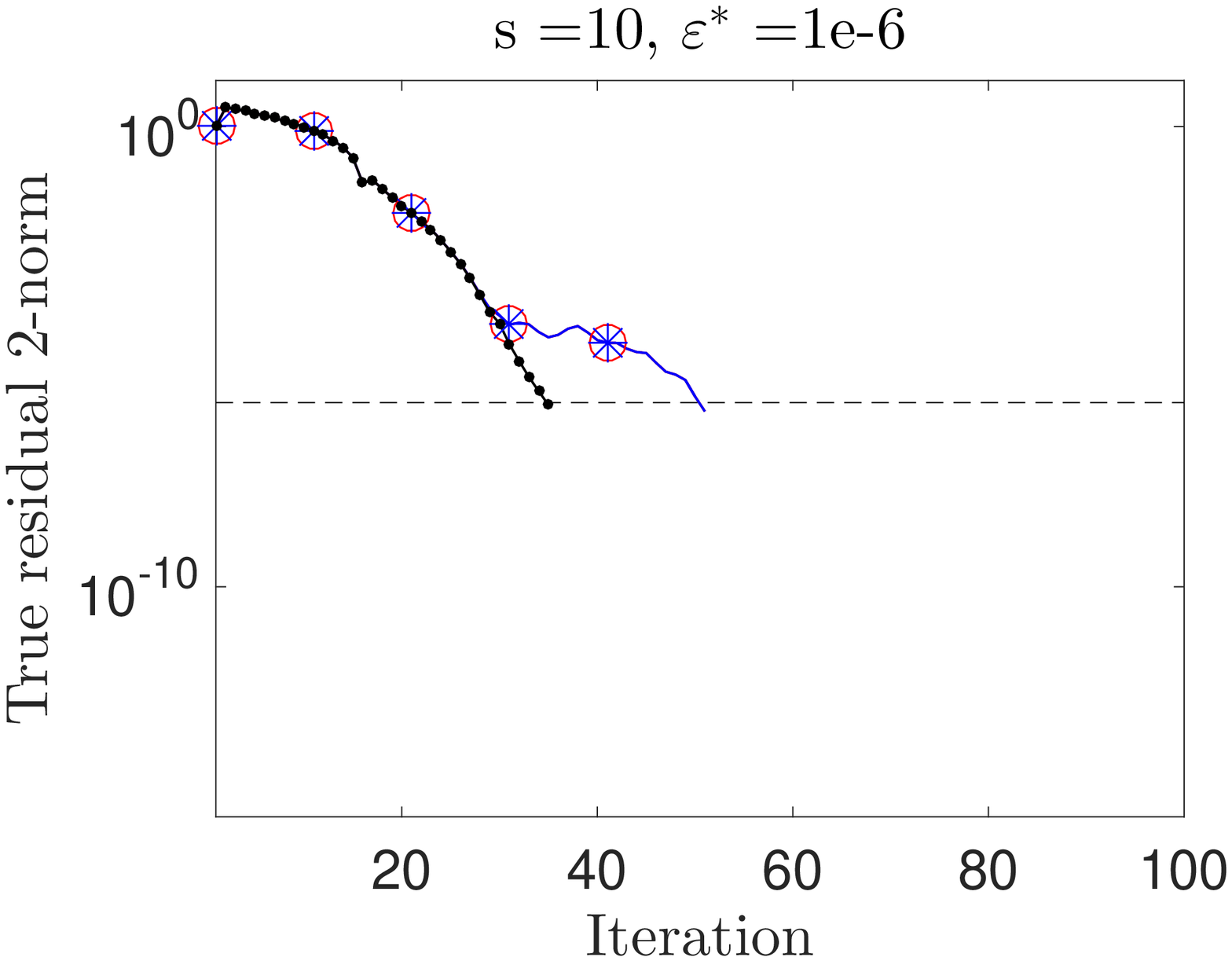}\\
	\label{fig:gr3030}
	
	\vspace{.6cm}
	\captionof{table}{Number of global synchronizations (outer loop iterations) in tests for matrix gr\_30\_30 (corresponding to Figure~\ref{fig:gr3030}). Dashes (-) in the table indicate that the method failed to converge to the desired level $\veps^*$.}
	\footnotesize
	\begin{tabular}{lr|c|c|c|}
\cline{3-5}
                                                         & \multicolumn{1}{l|}{} & fixed $s$-step CG & adaptive $s$-step CG & classical CG        \\ \hline
\multicolumn{1}{|l|}{\multirow{3}{*}{$\veps^*=$3.4e-14}} & $s=4$                 & 16          & 17               & \multirow{3}{*}{52} \\ \cline{2-4}
\multicolumn{1}{|l|}{}                                   & $s=8$                 & -           & 14               &                     \\ \cline{2-4}
\multicolumn{1}{|l|}{}                                   & $s=10$                & -           & 14               &                     \\ \hline
\multicolumn{1}{|l|}{\multirow{3}{*}{$\veps^*=$1e-6}}    & $s=4$                 & 9           & 9                & \multirow{3}{*}{34} \\ \cline{2-4}
\multicolumn{1}{|l|}{}                                   & $s=8$                 & 5           & 5                &                     \\ \cline{2-4}
\multicolumn{1}{|l|}{}                                   & $s=10$                & 5           & 5                &                     \\ \hline
\end{tabular}
\label{tab:gr3030}

\end{figure}

\begin{figure}[ht]
	\centering
			\captionof{figure}{Convergence of the 2-norm of the true residual for the matrix mesh3e1. Plots on the left use $\veps^*=$1.0e-14 and 
	plots on the right use $\veps^*=$1e-6. Plots are shown for three $s$ (or $s_{\text{max}}$) values: $4$ (top), $8$ (middle), and $10$ (bottom).  
	Iterations where communication takes place are noted with markers (red circles for $s$-step CG, blue stars for variable $s$-step CG, and black dots for classical CG). The horizontal dashed line shows the requested accuracy $\veps^*$. }
	\vspace{.2cm}
	\includegraphics[width=2.15in]{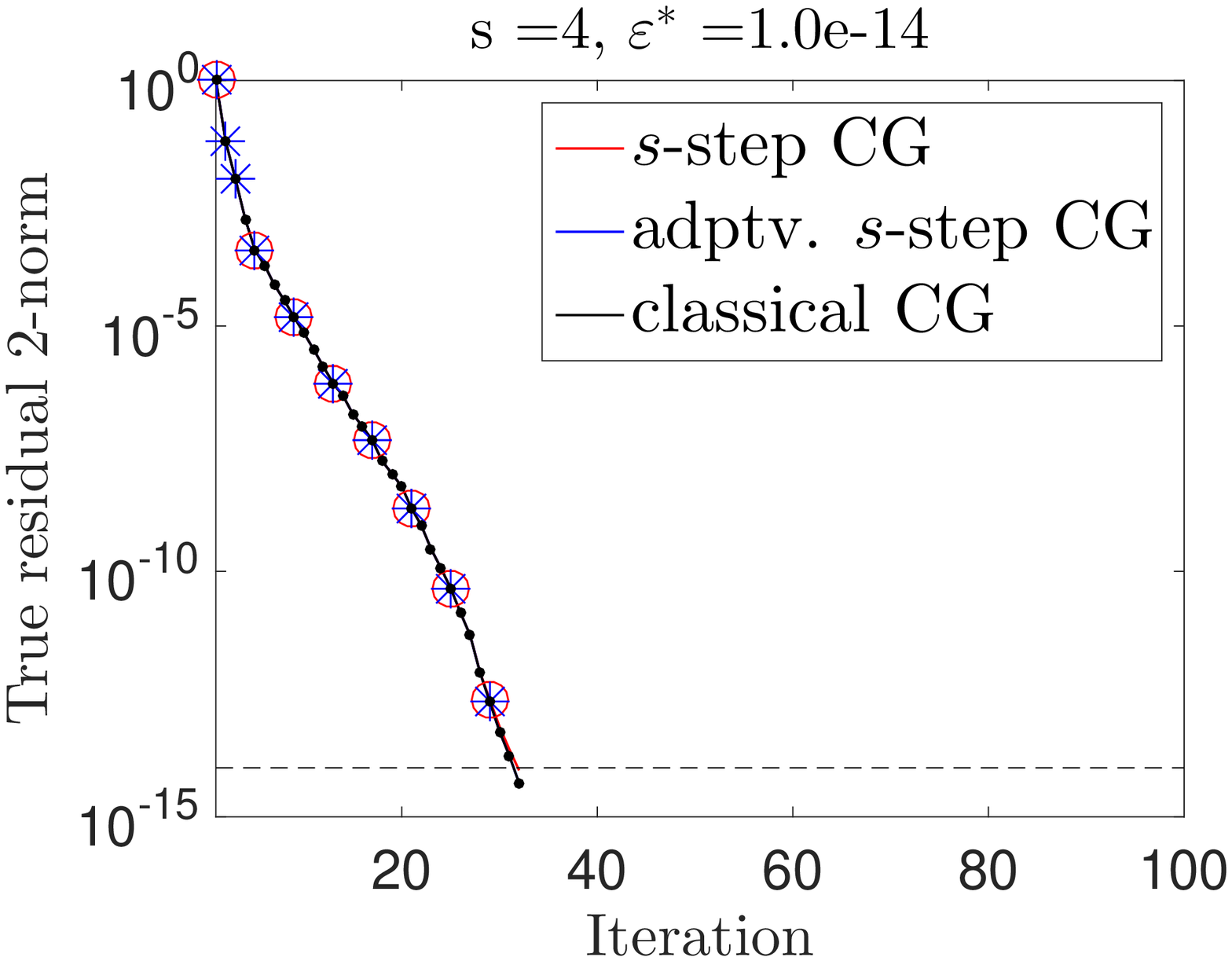}
	\hspace{5mm}
		\includegraphics[width=2.15in]{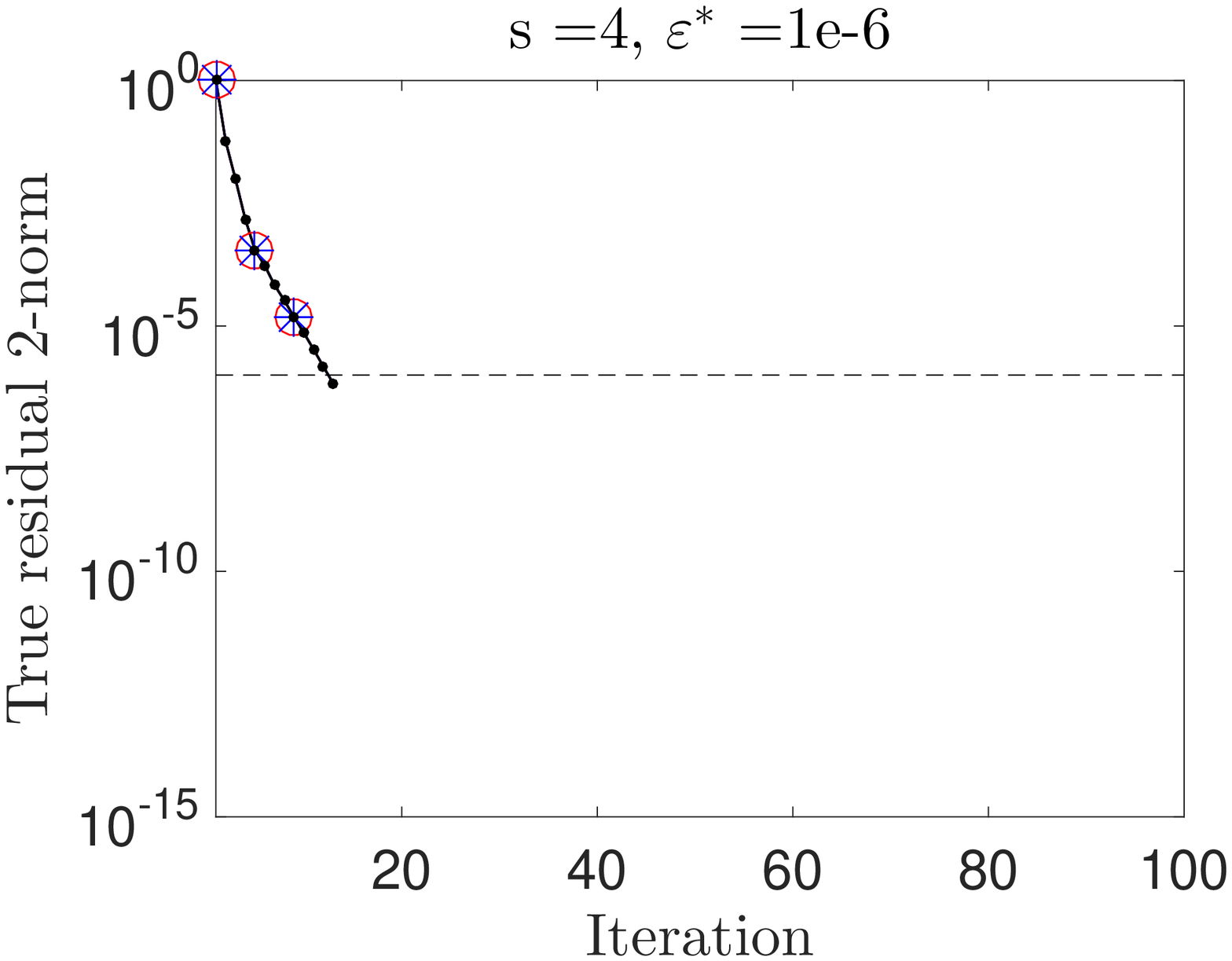}\\
		\vspace{1mm}
	\includegraphics[width=2.15in]{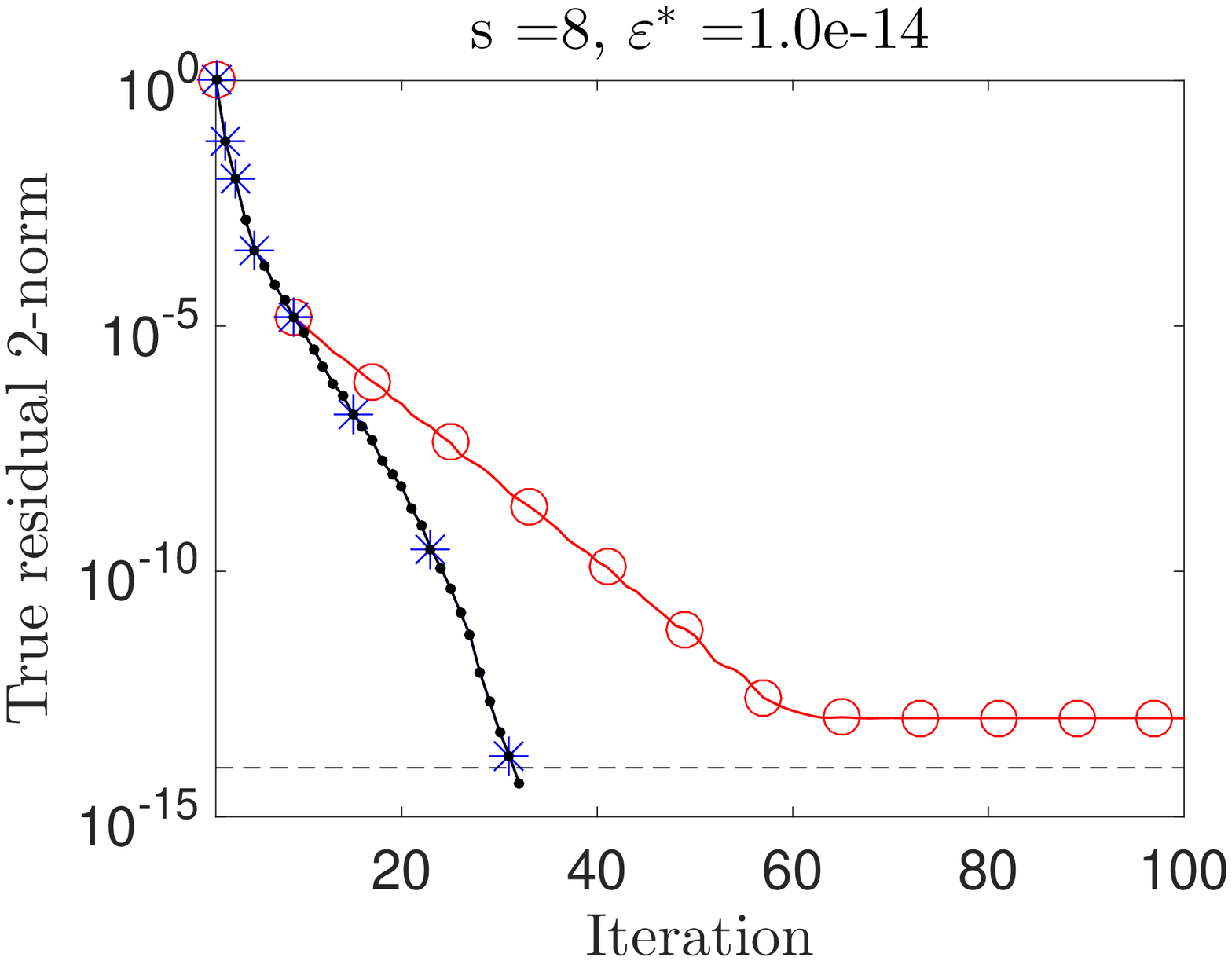}
	\hspace{5mm}
		\includegraphics[width=2.15in]{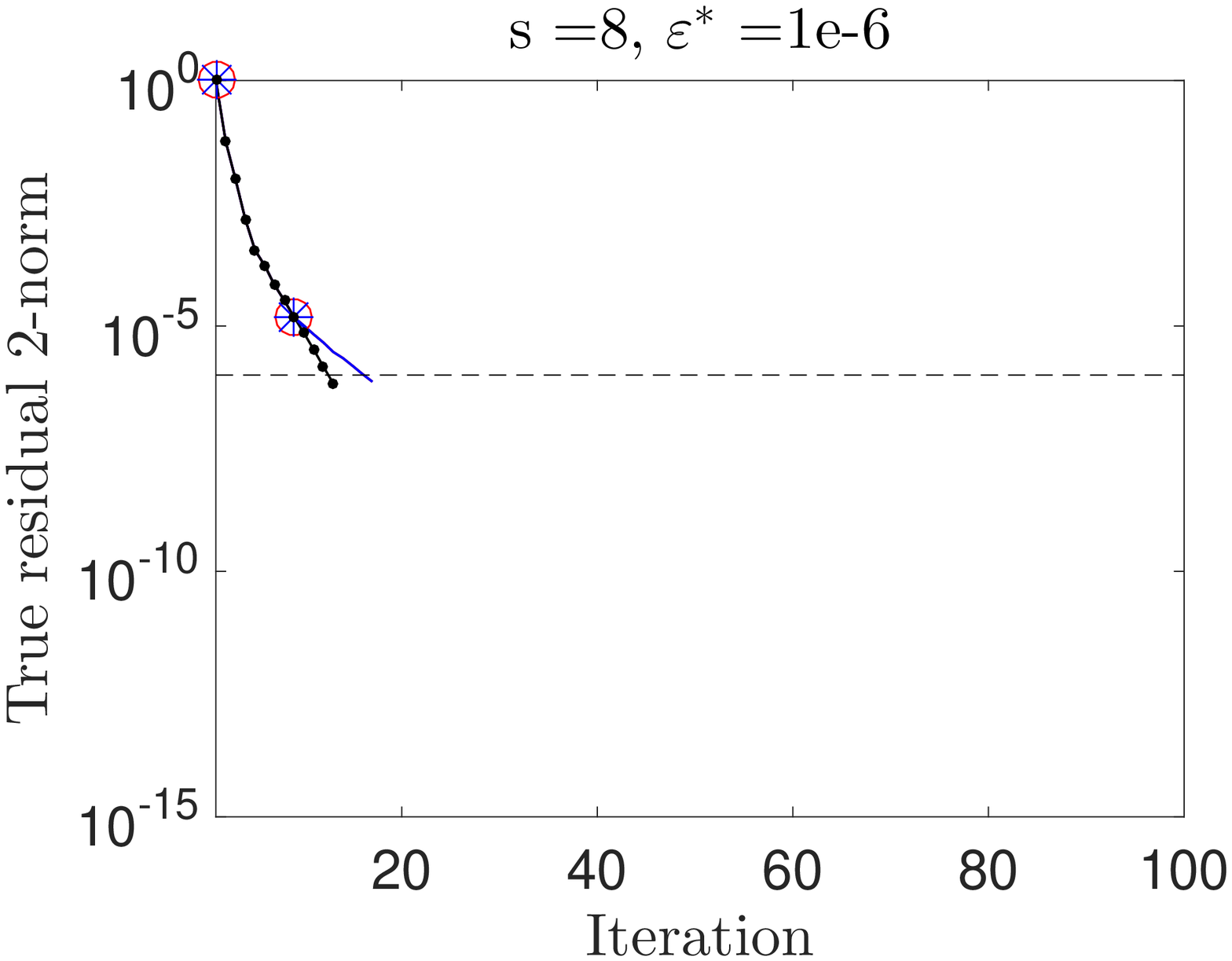}\\	
		\vspace{1mm}
			\includegraphics[width=2.15in]{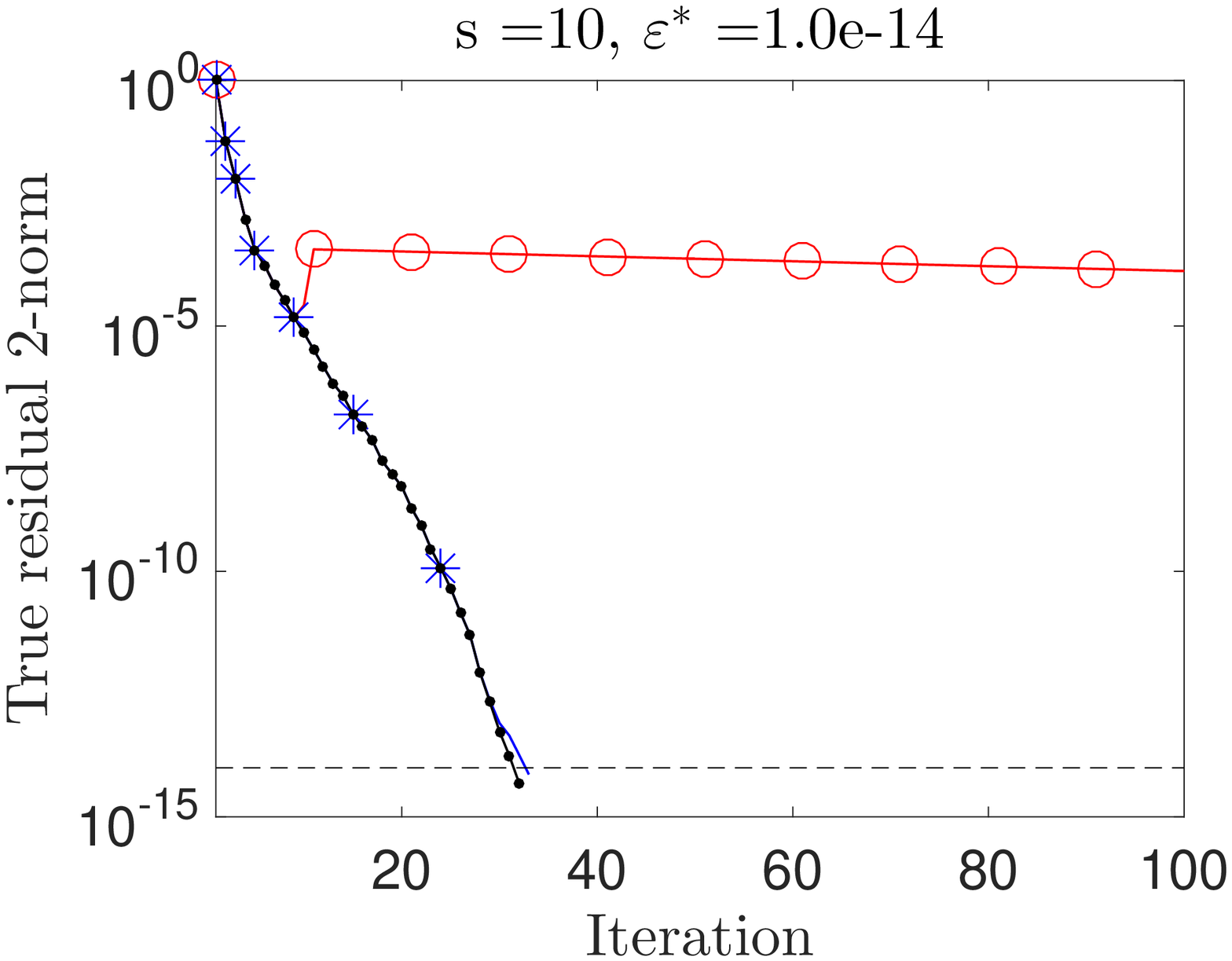}
	\hspace{5mm}
		\includegraphics[width=2.15in]{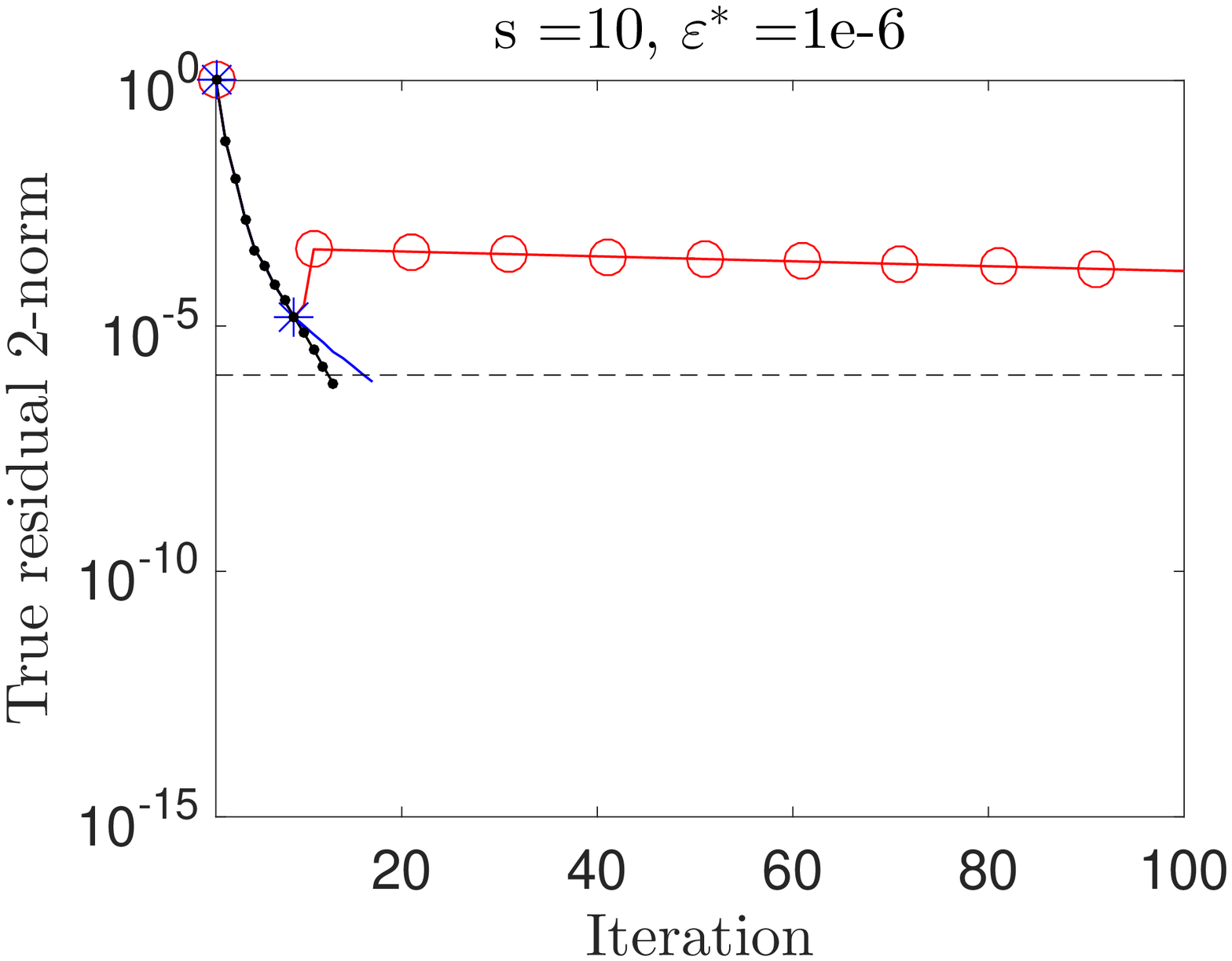}\\
	\label{fig:mesh3e1}
	
	\vspace{.6cm}
	\footnotesize
\captionof{table}{Number of global synchronizations (outer loop iterations) in tests for matrix mesh3e1 (corresponding to Figure~\ref{fig:mesh3e1}). Dashes (-) in the table indicate that the method failed to converge to the desired level $\veps^*$.}
\begin{tabular}{lr|c|c|c|}
\cline{3-5}
                                                       & \multicolumn{1}{l|}{} & fixed $s$-step CG & adaptive $s$-step CG & classical CG        \\ \hline
\multicolumn{1}{|l|}{\multirow{3}{*}{$\veps^*=$1e-14}} & $s=4$                 & 8           & 10               & \multirow{3}{*}{31} \\ \cline{2-4}
\multicolumn{1}{|l|}{}                                 & $s=8$                 & -           & 8                &                     \\ \cline{2-4}
\multicolumn{1}{|l|}{}                                 & $s=10$                & -           & 7                &                     \\ \hline
\multicolumn{1}{|l|}{\multirow{3}{*}{$\veps^*=$1e-6}}  & $s=4$                 & 3           & 3                & \multirow{3}{*}{12} \\ \cline{2-4}
\multicolumn{1}{|l|}{}                                 & $s=8$                 & 2           & 2                &                     \\ \cline{2-4}
\multicolumn{1}{|l|}{}                                 & $s=10$                & 99          & 2                &                     \\ \hline
\end{tabular}
\label{tab:mesh3e1}

\end{figure}

\begin{figure}[ht]
	\centering
			\captionof{figure}{Convergence of the 2-norm of the true residual for the matrix nos6. Plots on the left use $\veps^*=$5.5e-10 and 
	plots on the right use $\veps^*=$1e-6. Plots are shown for three $s$ (or $s_{\text{max}}$) values: $4$ (top), $8$ (middle), and $10$ (bottom).  
	Iterations where communication takes place are noted with markers (red circles for $s$-step CG, blue stars for variable $s$-step CG, and black dots for classical CG). The horizontal dashed line shows the requested accuracy $\veps^*$. }
	\vspace{.2cm}
	\includegraphics[width=2.15in]{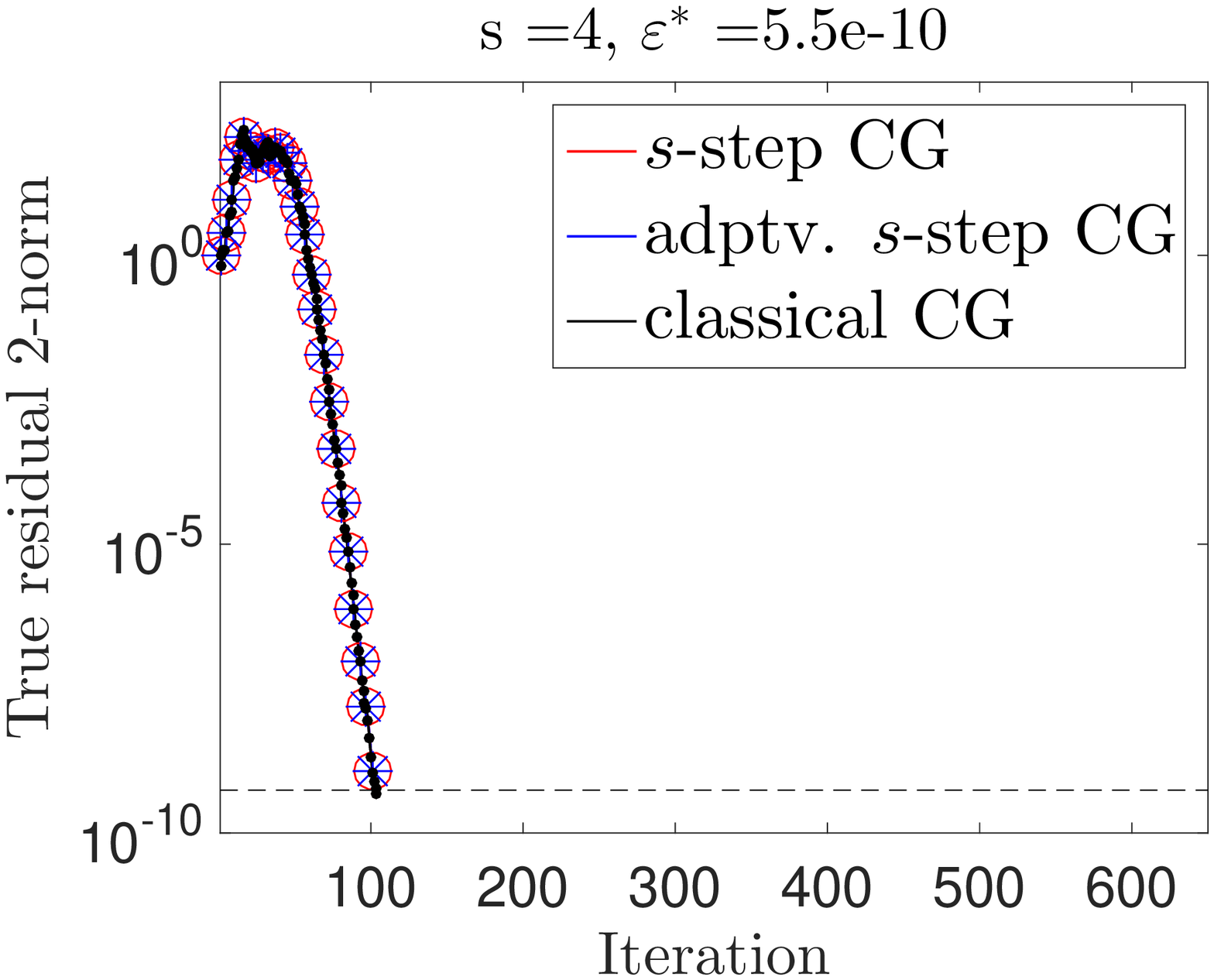}
	\hspace{5mm}
		\includegraphics[width=2.15in]{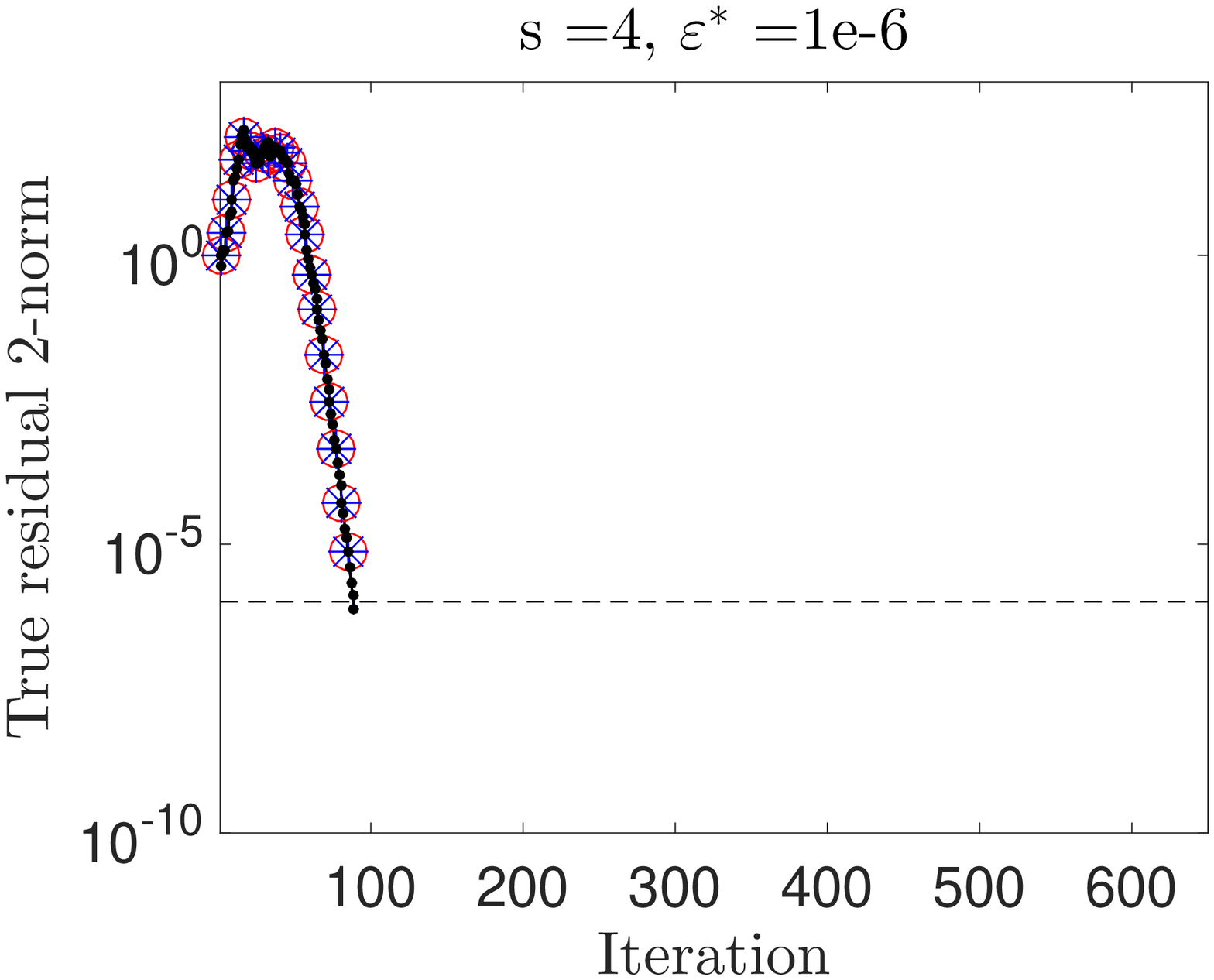}\\
		\vspace{1mm}
	\includegraphics[width=2.15in]{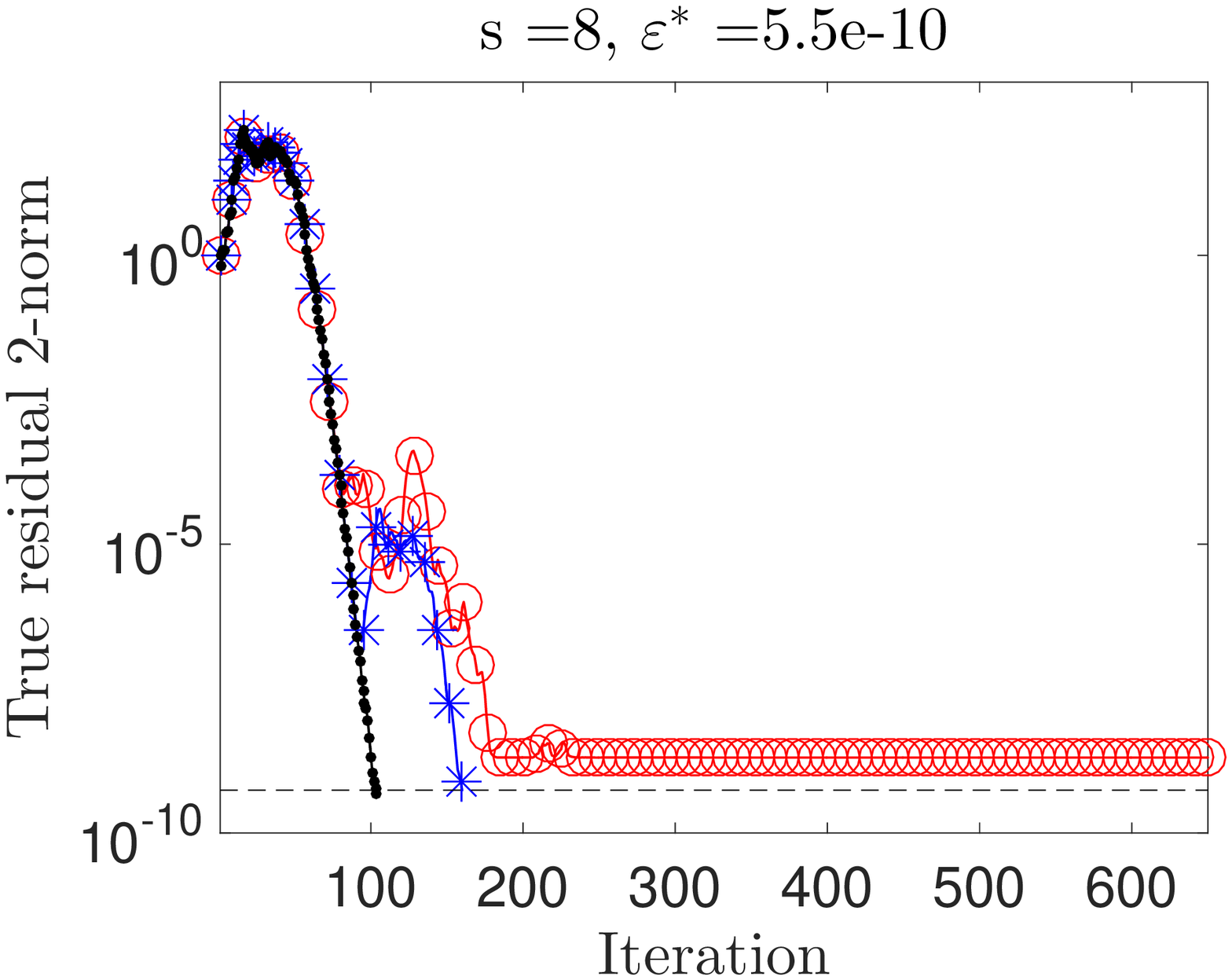}
	\hspace{5mm}
		\includegraphics[width=2.15in]{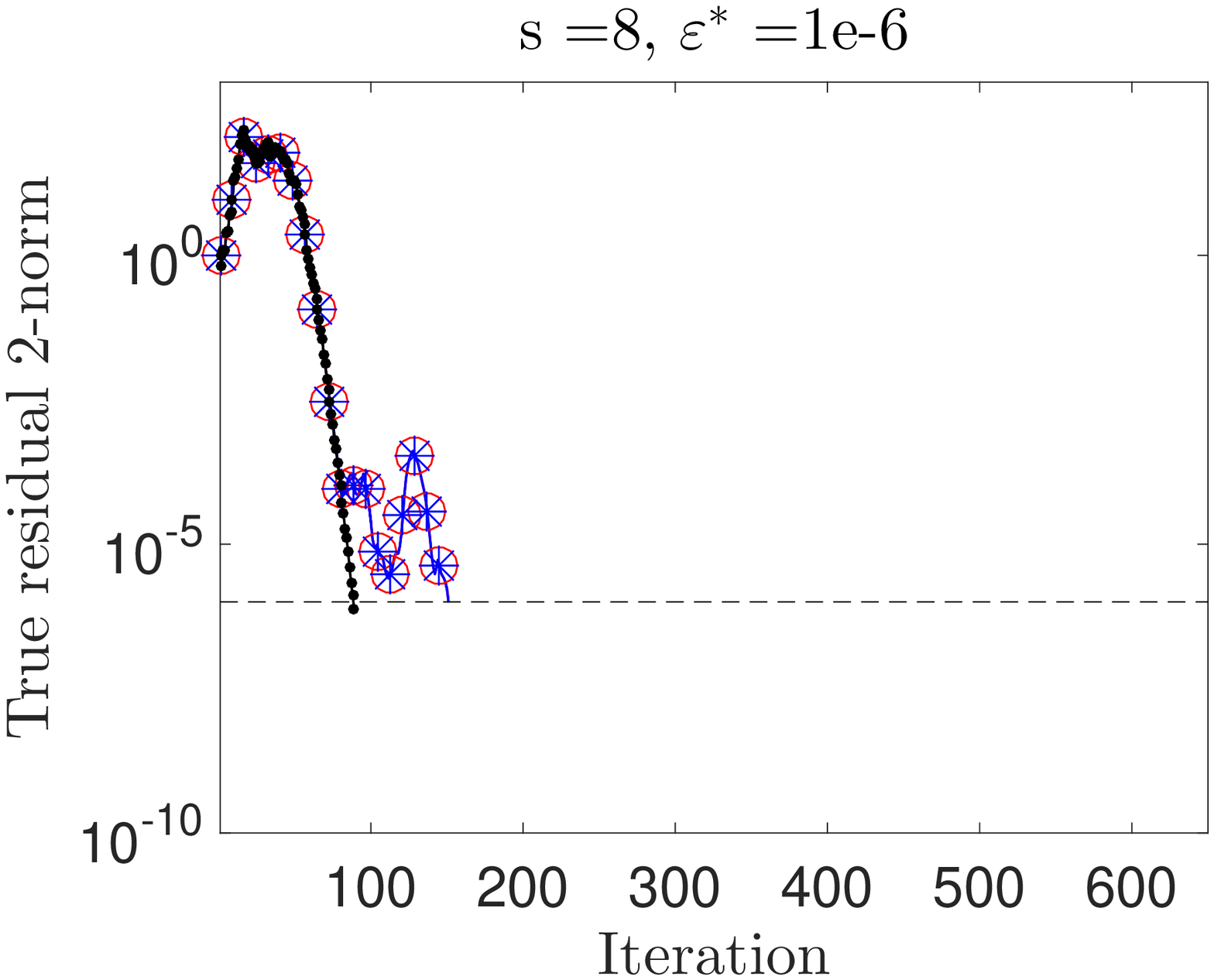}\\	
		\vspace{1mm}
			\includegraphics[width=2.15in]{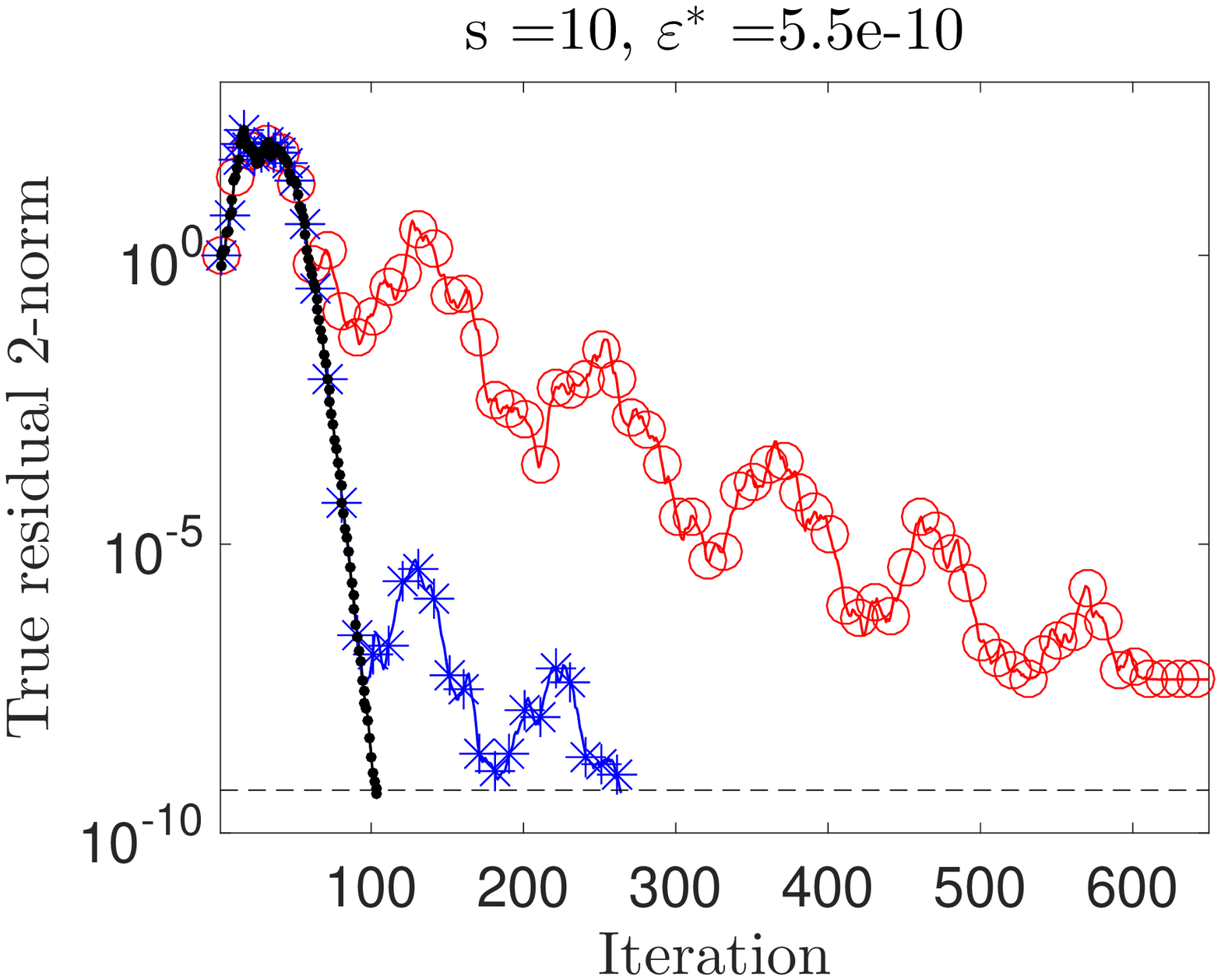}
	\hspace{5mm}
		\includegraphics[width=2.15in]{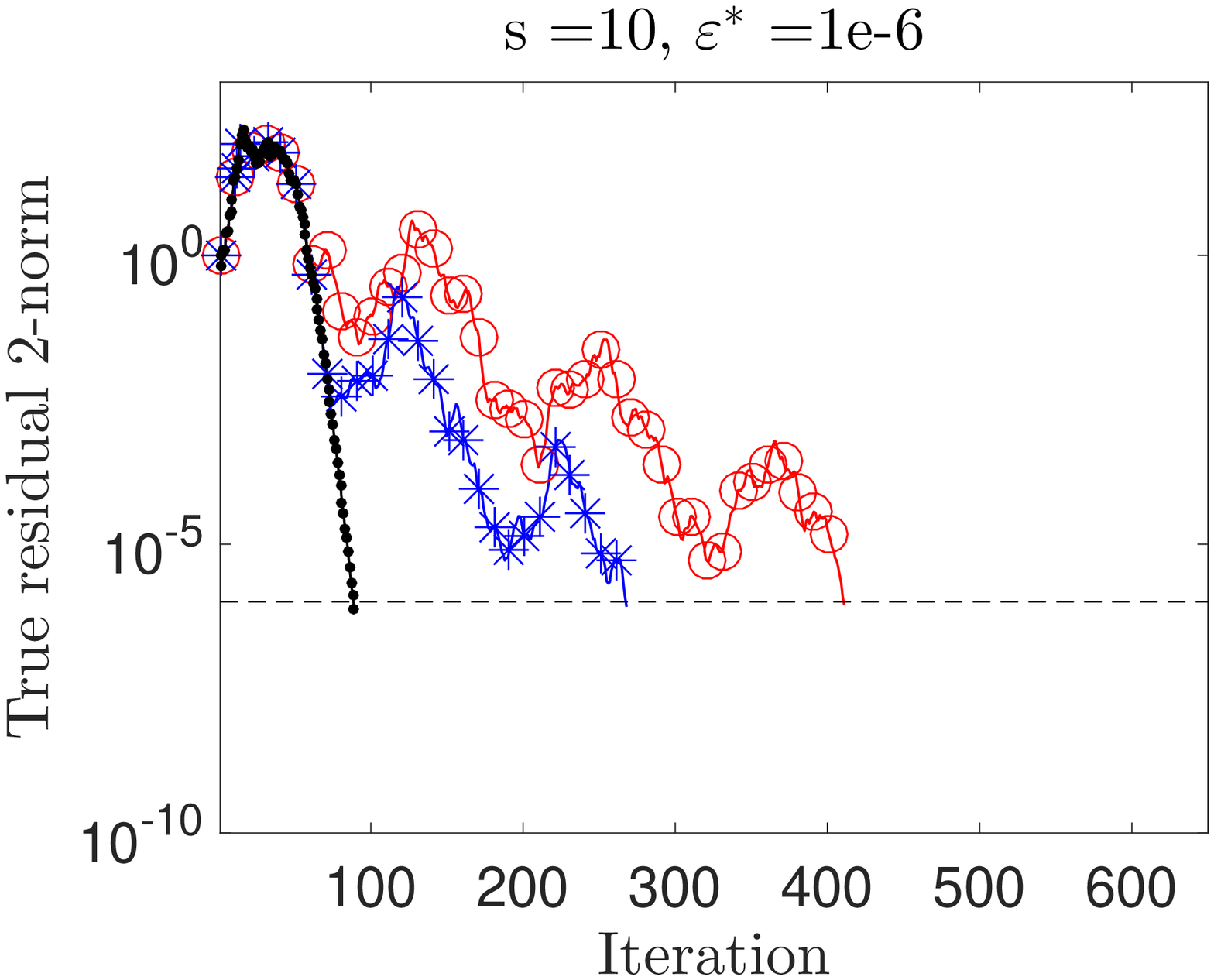}\\
	\label{fig:nos6}
	
	\vspace{.6cm}
	\footnotesize
\captionof{table}{Number of global synchronizations (outer loop iterations) in tests for matrix nos6 (corresponding to Figure~\ref{fig:nos6}). Dashes (-) in the table indicate that the method failed to converge to the desired level $\veps^*$.}
\begin{tabular}{lr|c|c|c|}
\cline{3-5}
                                                       & \multicolumn{1}{l|}{} & fixed $s$-step CG & adaptive $s$-step CG & classical CG         \\ \hline
\multicolumn{1}{|l|}{\multirow{3}{*}{$\veps^*=$5.5e-10}} & $s=4$                 & 26          & 26               & \multirow{3}{*}{103} \\ \cline{2-4}
\multicolumn{1}{|l|}{}                                 & $s=8$                 & -           & 29               &                      \\ \cline{2-4}
\multicolumn{1}{|l|}{}                                 & $s=10$                & -           & 36               &                      \\ \hline
\multicolumn{1}{|l|}{\multirow{3}{*}{$\veps^*=$1e-6}}  & $s=4$                 & 22          & 22               & \multirow{3}{*}{88}  \\ \cline{2-4}
\multicolumn{1}{|l|}{}                                 & $s=8$                 & 19          & 19               &                      \\ \cline{2-4}
\multicolumn{1}{|l|}{}                                 & $s=10$                & 41          & 29               &                      \\ \hline
\end{tabular}
\label{tab:nos6}
\end{figure}

\begin{figure}[ht]
	\centering
			\captionof{figure}{Convergence of the 2-norm of the true residual for the matrix bcsstk09. Plots on the left use $\veps^*=$1.6e-12 and 
	plots on the right use $\veps^*=$1e-6. Plots are shown for three $s$ (or $s_{\text{max}}$) values: $4$ (top), $8$ (middle), and $10$ (bottom).  
	Iterations where communication takes place are noted with markers (red circles for $s$-step CG, blue stars for variable $s$-step CG, and black dots for classical CG). The horizontal dashed line shows the requested accuracy $\veps^*$.}
	\vspace{.2cm}
	\includegraphics[width=2.15in]{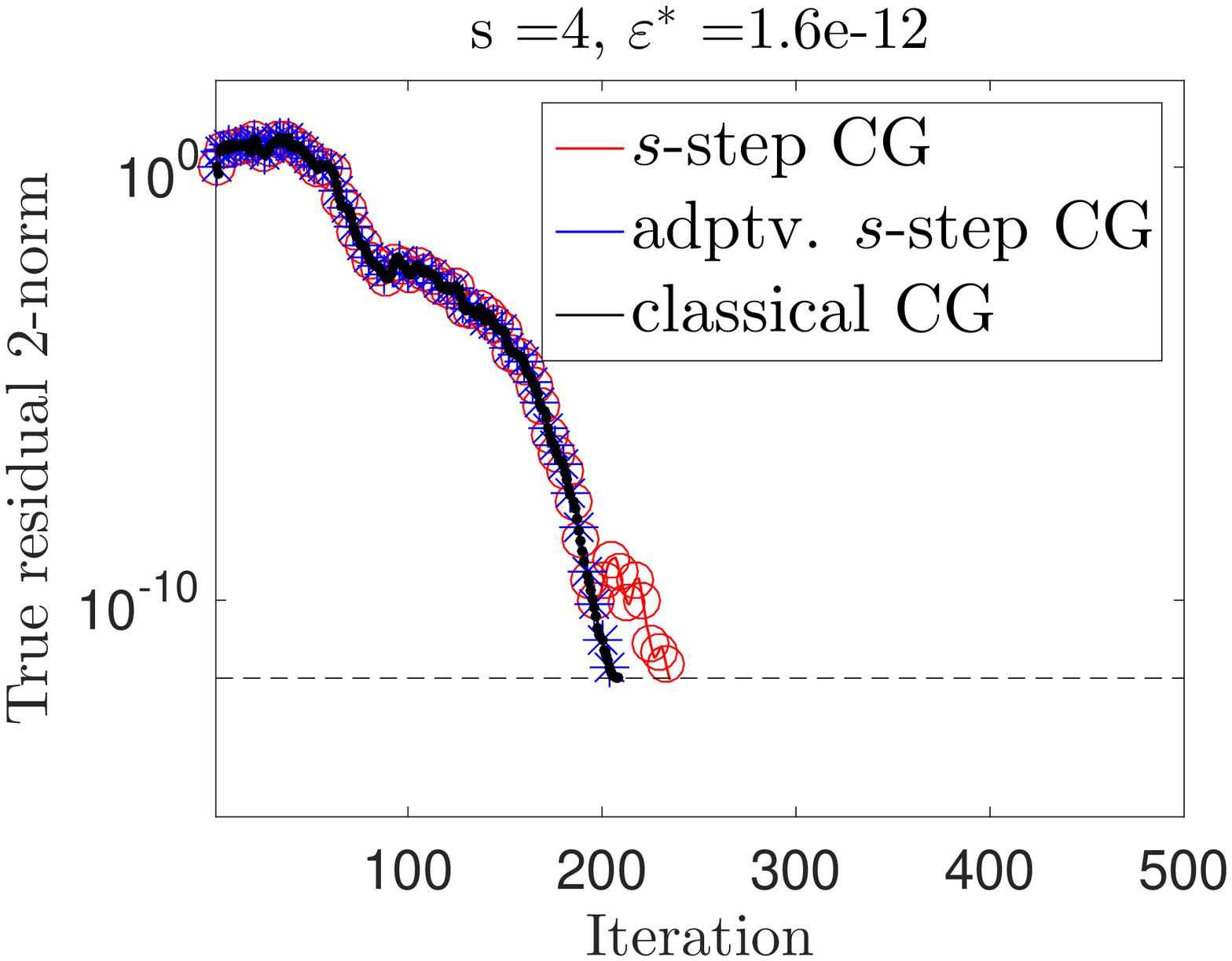}
	\hspace{5mm}
		\includegraphics[width=2.15in]{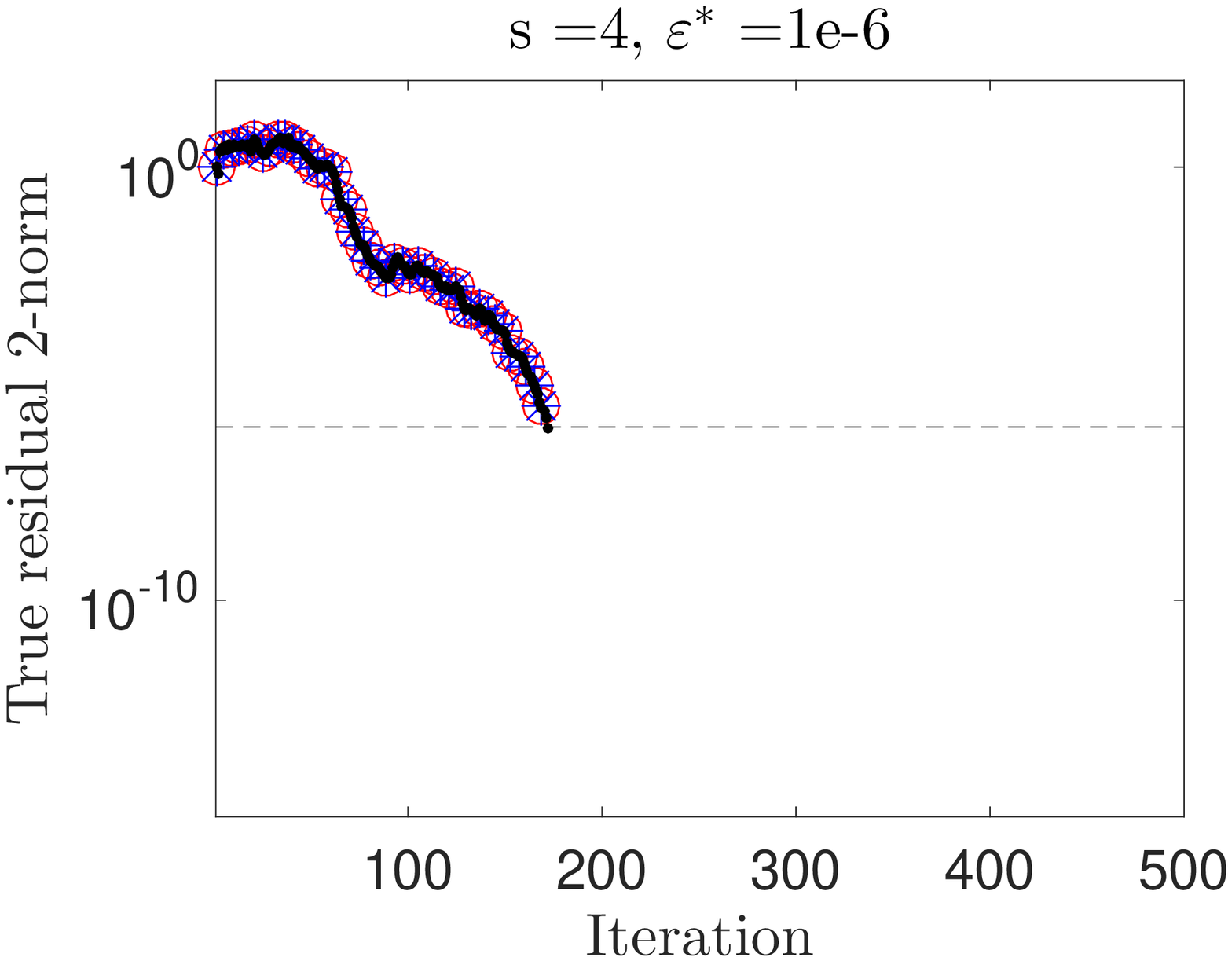}\\
		\vspace{1mm}
	\includegraphics[width=2.15in]{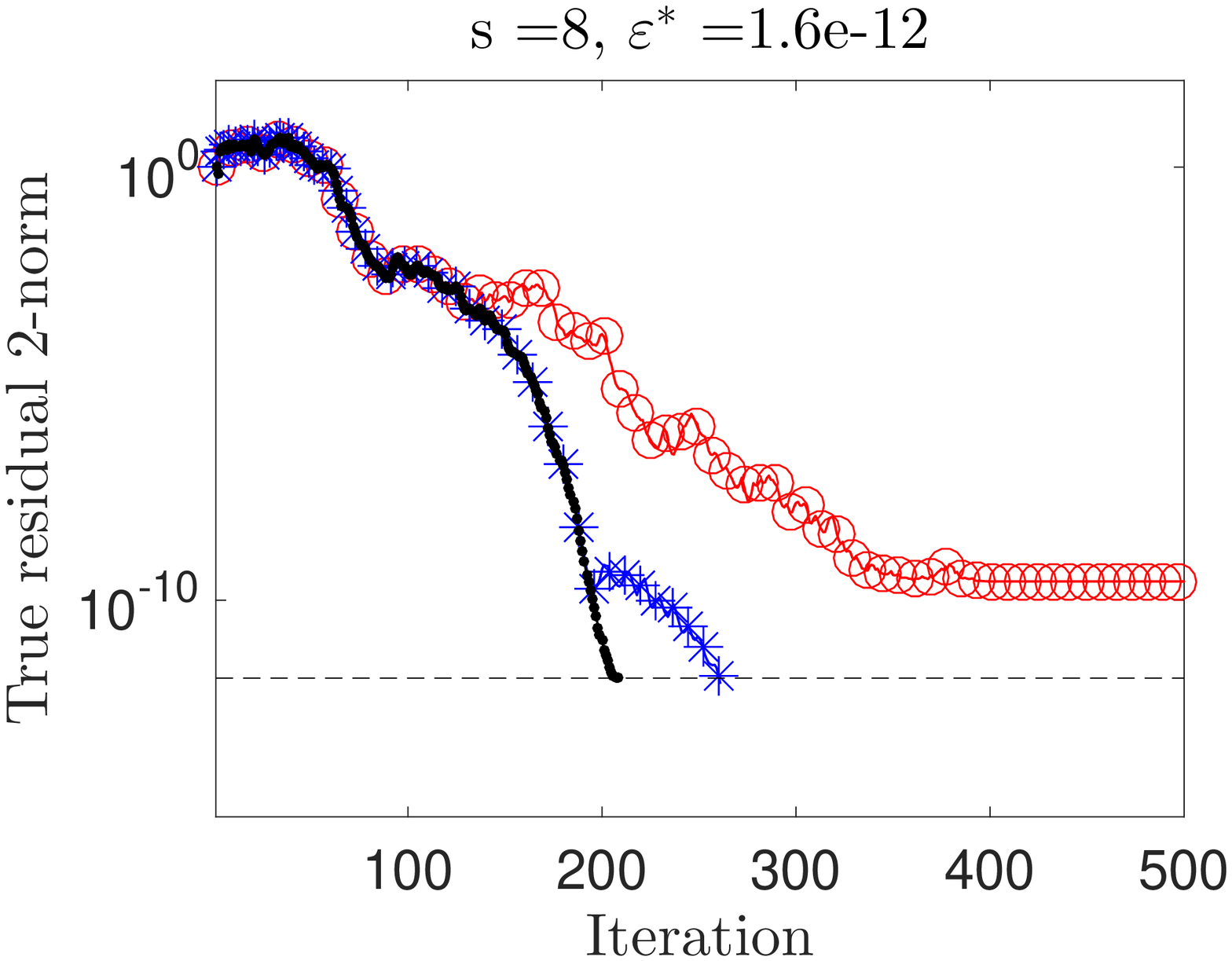}
	\hspace{5mm}
		\includegraphics[width=2.15in]{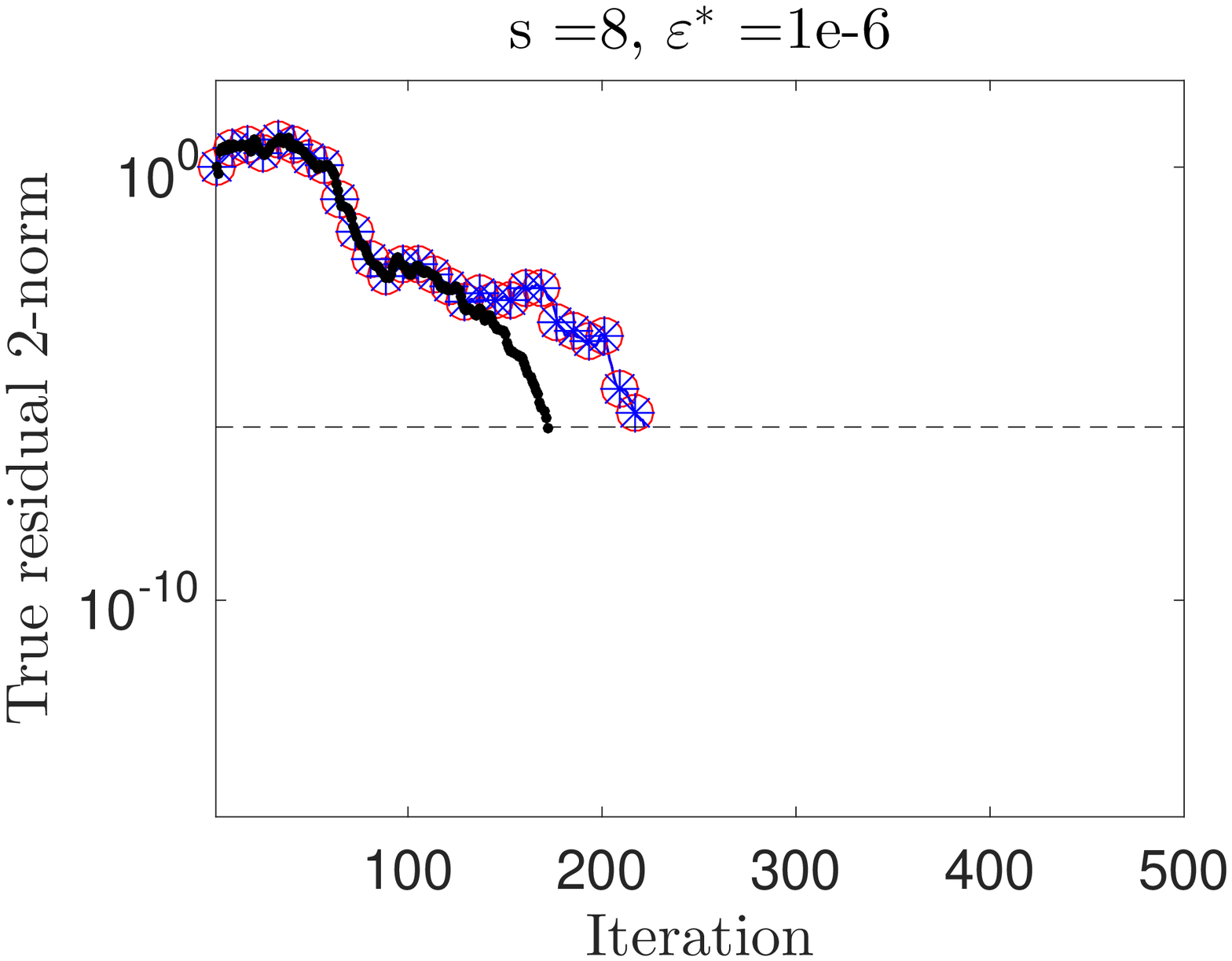}\\
			\vspace{1mm}
			\includegraphics[width=2.15in]{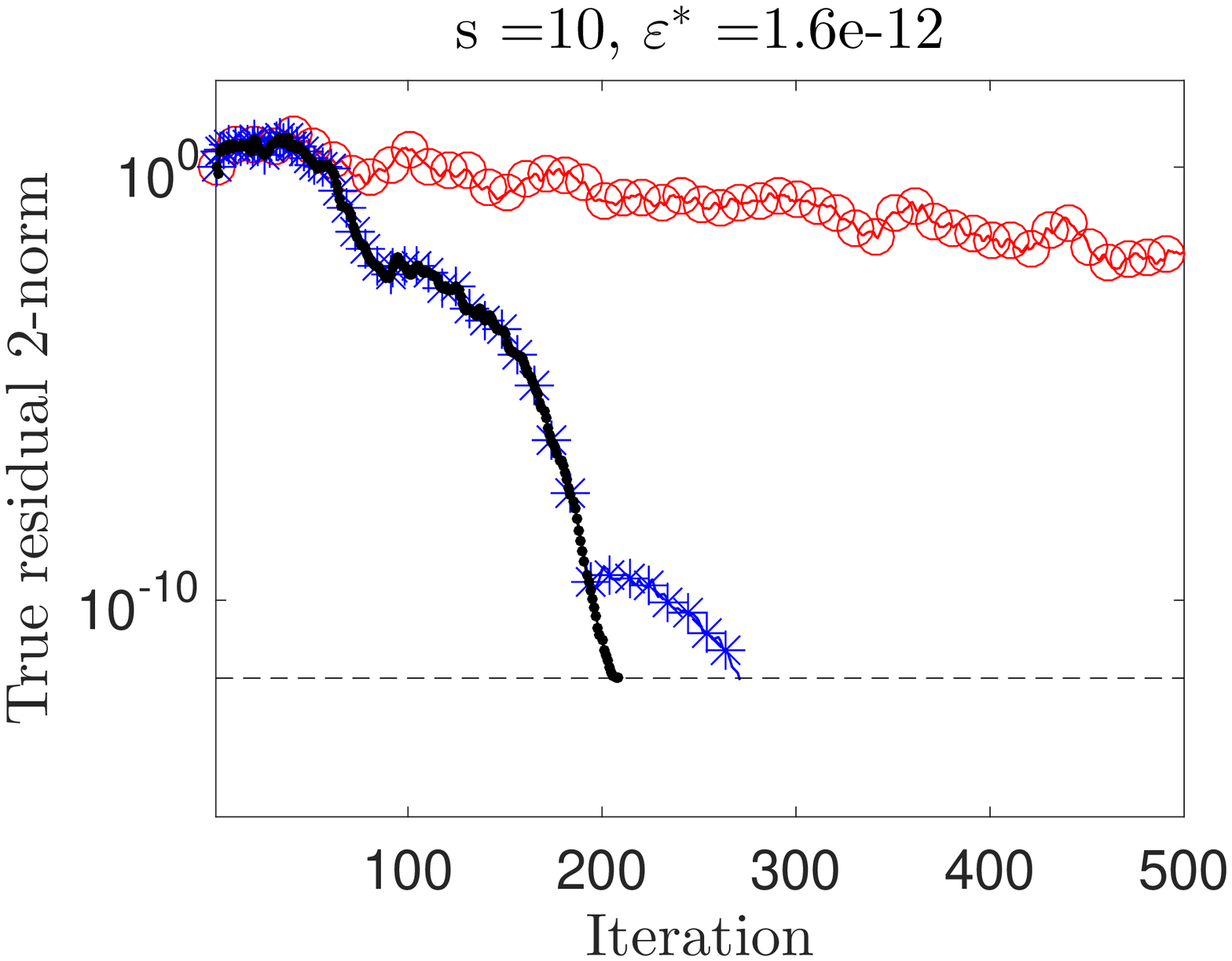}
	\hspace{5mm}
		\includegraphics[width=2.15in]{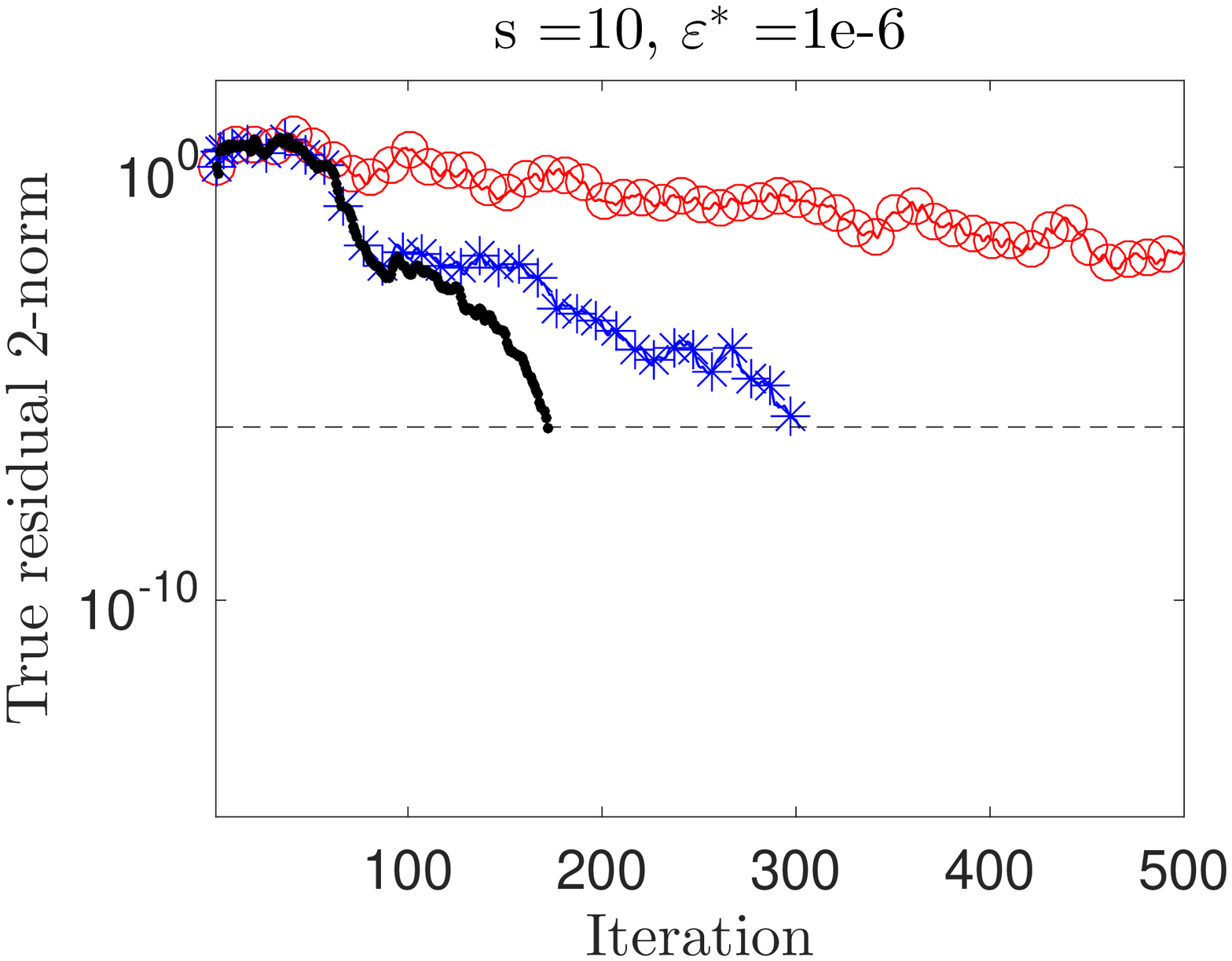}\\
	\label{fig:bcsstk09}
	
	\vspace{.6cm}
	\footnotesize
\captionof{table}{Number of global synchronizations (outer loop iterations) in tests for matrix bcsstk09 (corresponding to Figure~\ref{fig:bcsstk09}). Dashes (-) in the table indicate that the method failed to converge to the desired level $\veps^*$. Tests used $c_{k,j}=10$.}
\begin{tabular}{lr|c|c|c|}
\cline{3-5}
                                                         & \multicolumn{1}{l|}{} & fixed $s$-step CG & adaptive $s$-step CG & classical CG         \\ \hline
\multicolumn{1}{|l|}{\multirow{3}{*}{$\veps^*=$1.6e-12}} & $s=4$                 & 59          & 59               & \multirow{3}{*}{207} \\ \cline{2-4}
\multicolumn{1}{|l|}{}                                   & $s=8$                 & -           & 52               &                      \\ \cline{2-4}
\multicolumn{1}{|l|}{}                                   & $s=10$                & -           & 50               &                      \\ \hline
\multicolumn{1}{|l|}{\multirow{3}{*}{$\veps^*=$1e-6}}    & $s=4$                 & 43          & 43               & \multirow{3}{*}{171} \\ \cline{2-4}
\multicolumn{1}{|l|}{}                                   & $s=8$                 & 28          & 28               &                      \\ \cline{2-4}
\multicolumn{1}{|l|}{}                                   & $s=10$                & 122         & 33               &                      \\ \hline
\end{tabular}
\label{tab:bcsstk09}

\end{figure}

\begin{figure}[ht]
	\centering
			\captionof{figure}{Convergence of the 2-norm of the true residual for the matrix ex5. Plots on the left use $\veps^*=$1e-8 and 
	plots on the right use $\veps^*=$1e-6. Plots are shown for three $s$ (or $s_{\text{max}}$) values: $4$ (top), $8$ (middle), and $10$ (bottom).  
	Iterations where communication takes place are noted with markers (red circles for $s$-step CG, blue stars for variable $s$-step CG, and black dots for classical CG). The horizontal dashed line shows the requested accuracy $\veps^*$.}
	\vspace{.2cm}
	\includegraphics[width=2.15in]{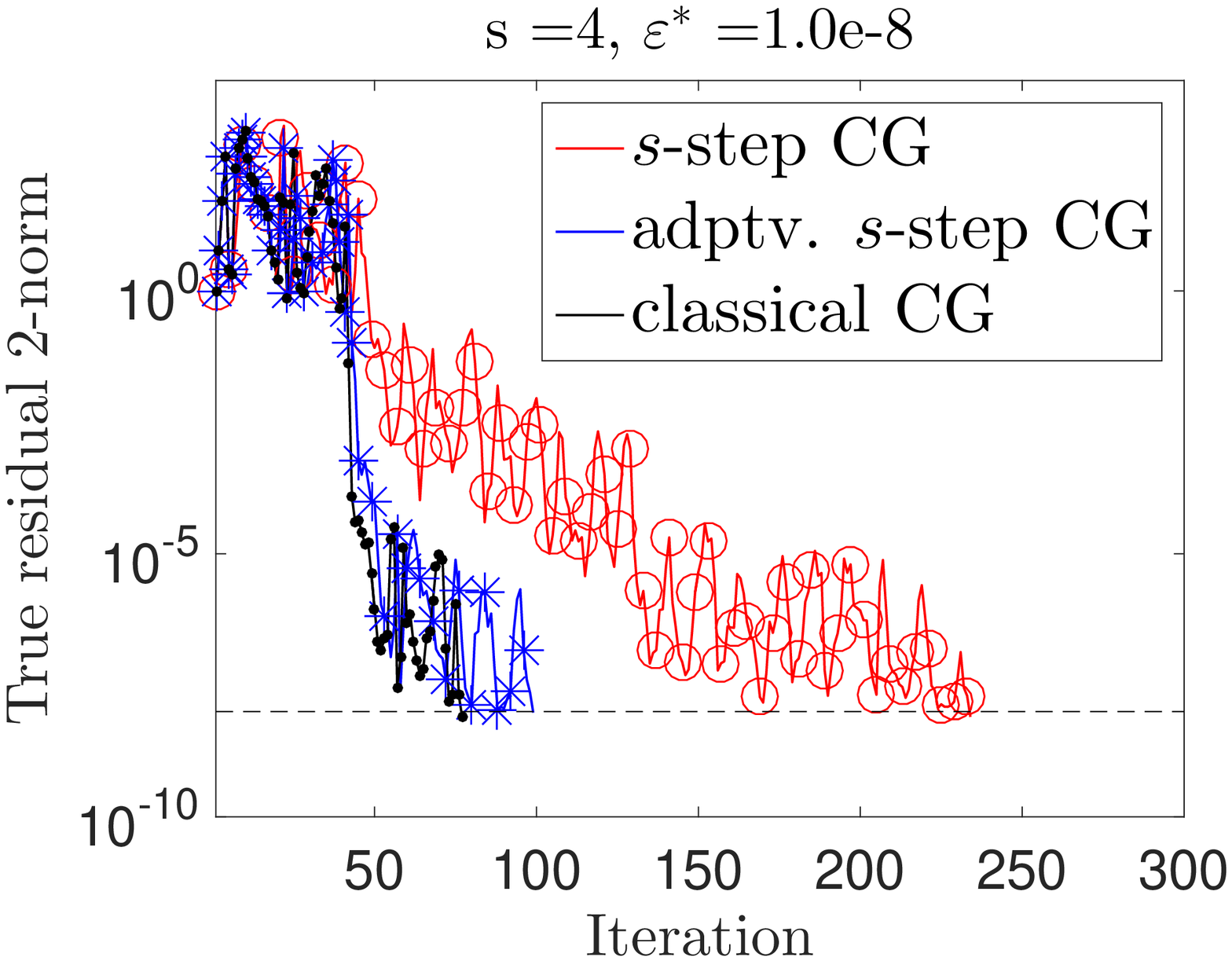}
	\hspace{5mm}
		\includegraphics[width=2.15in]{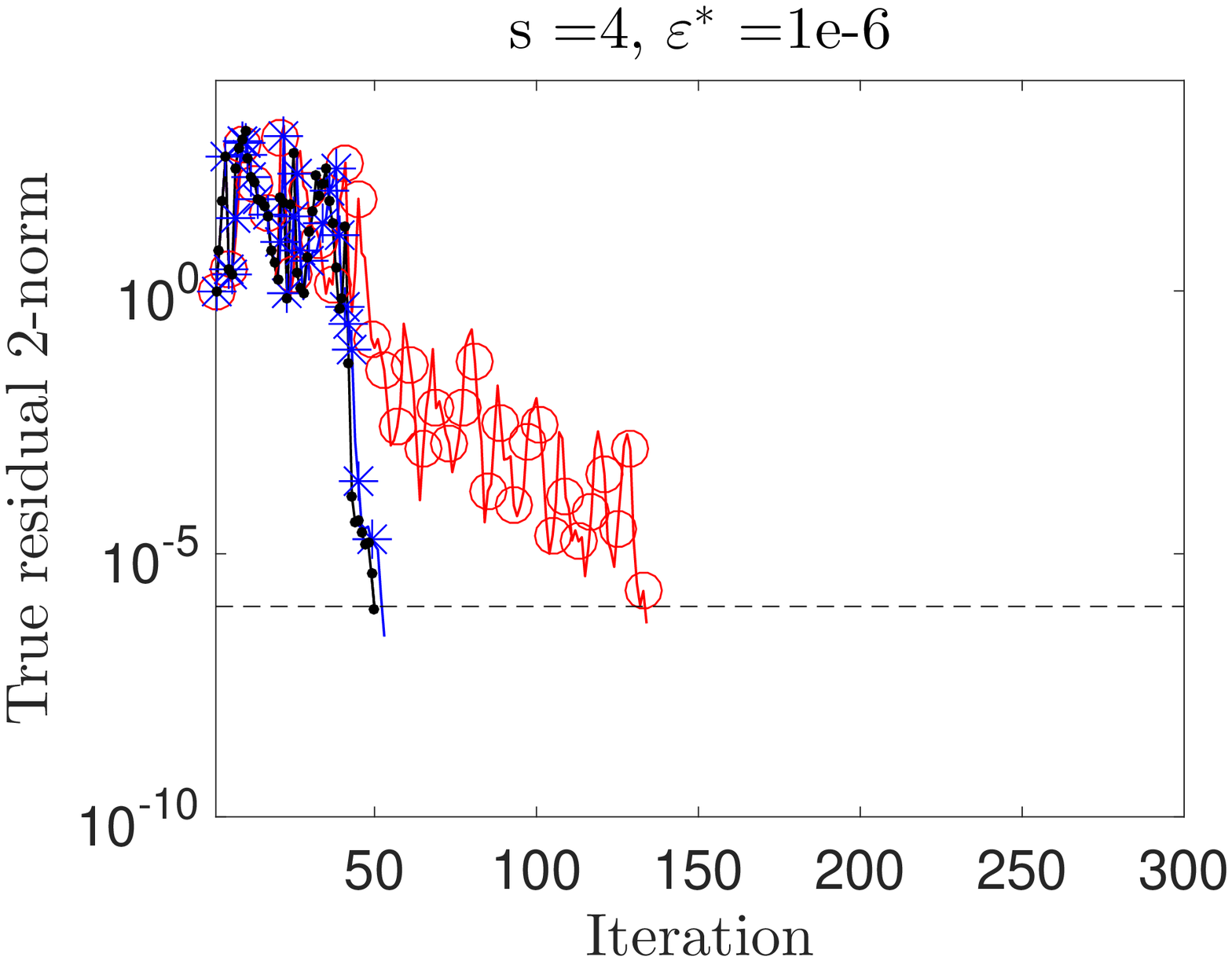}\\
		\vspace{1mm}
	\includegraphics[width=2.15in]{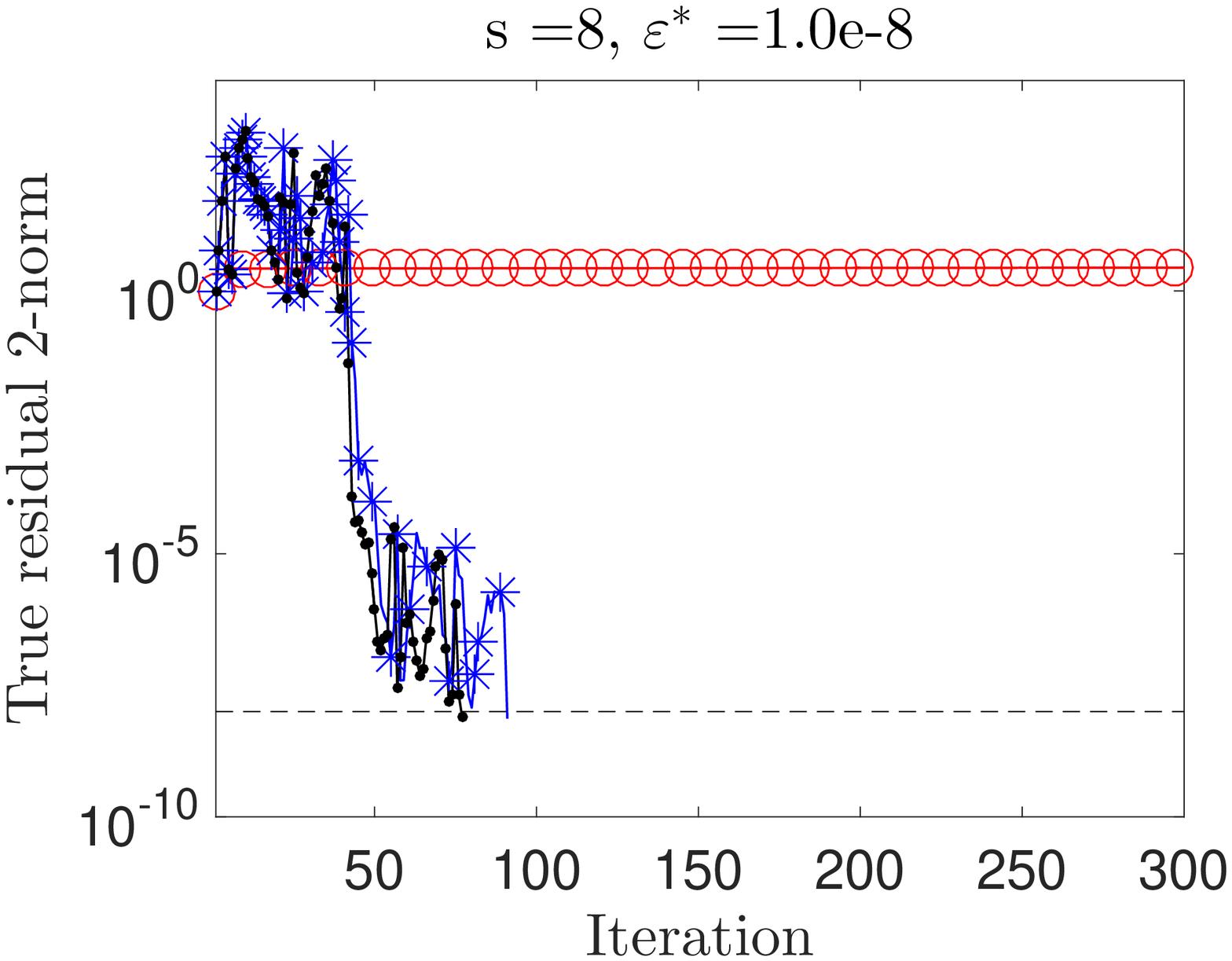}
	\hspace{5mm}
		\includegraphics[width=2.15in]{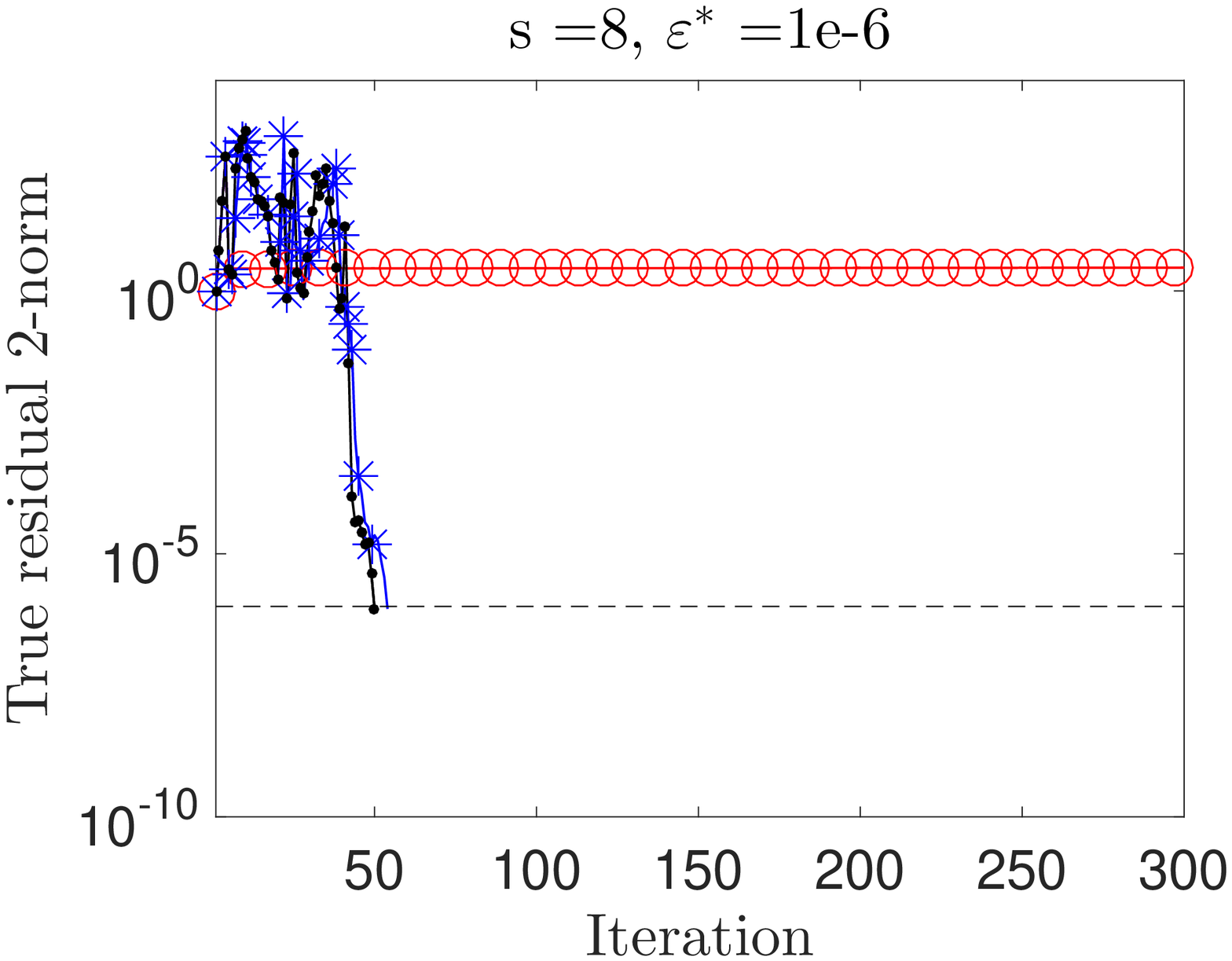}\\	
		\vspace{1mm}
			\includegraphics[width=2.15in]{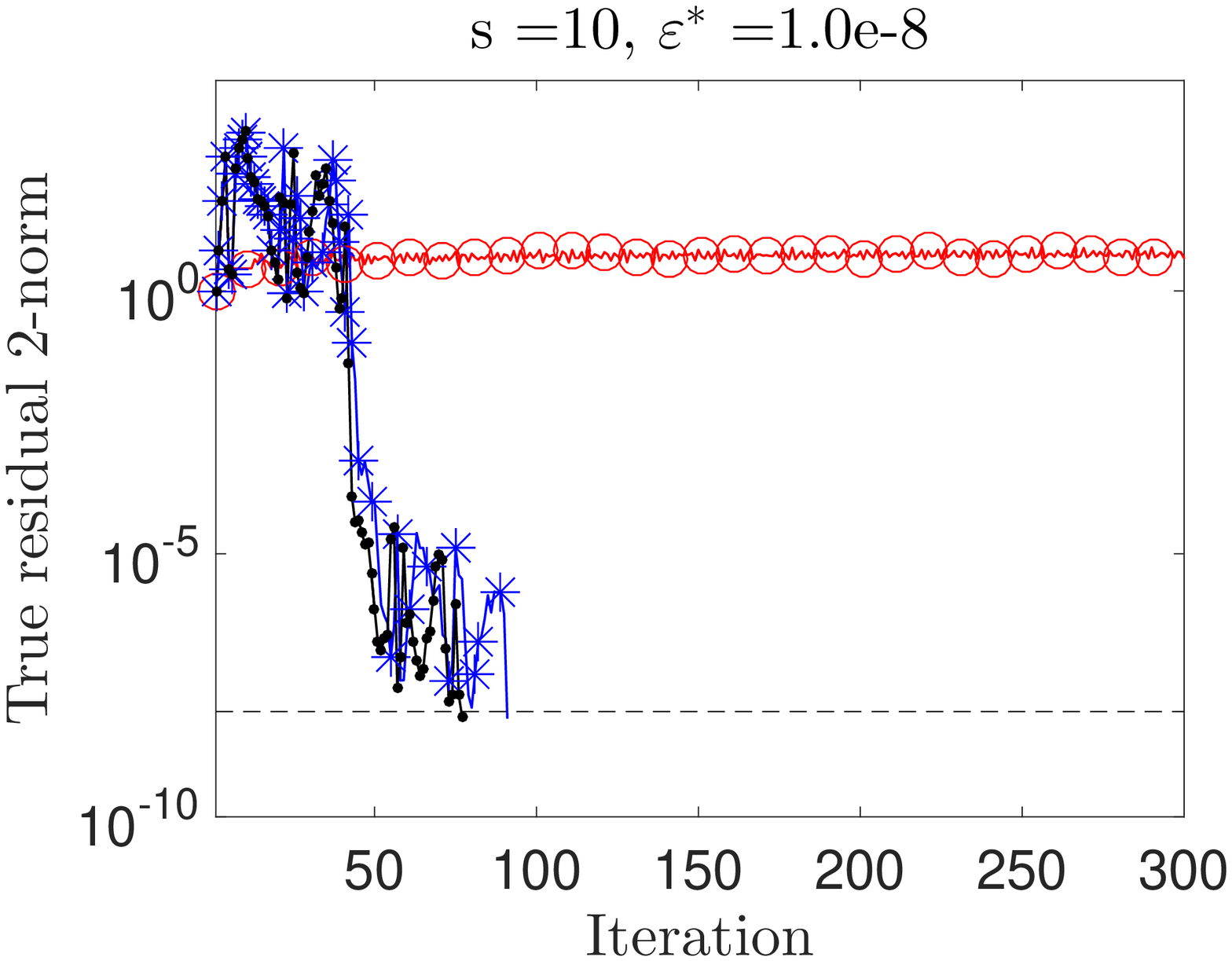}
	\hspace{5mm}
		\includegraphics[width=2.1in]{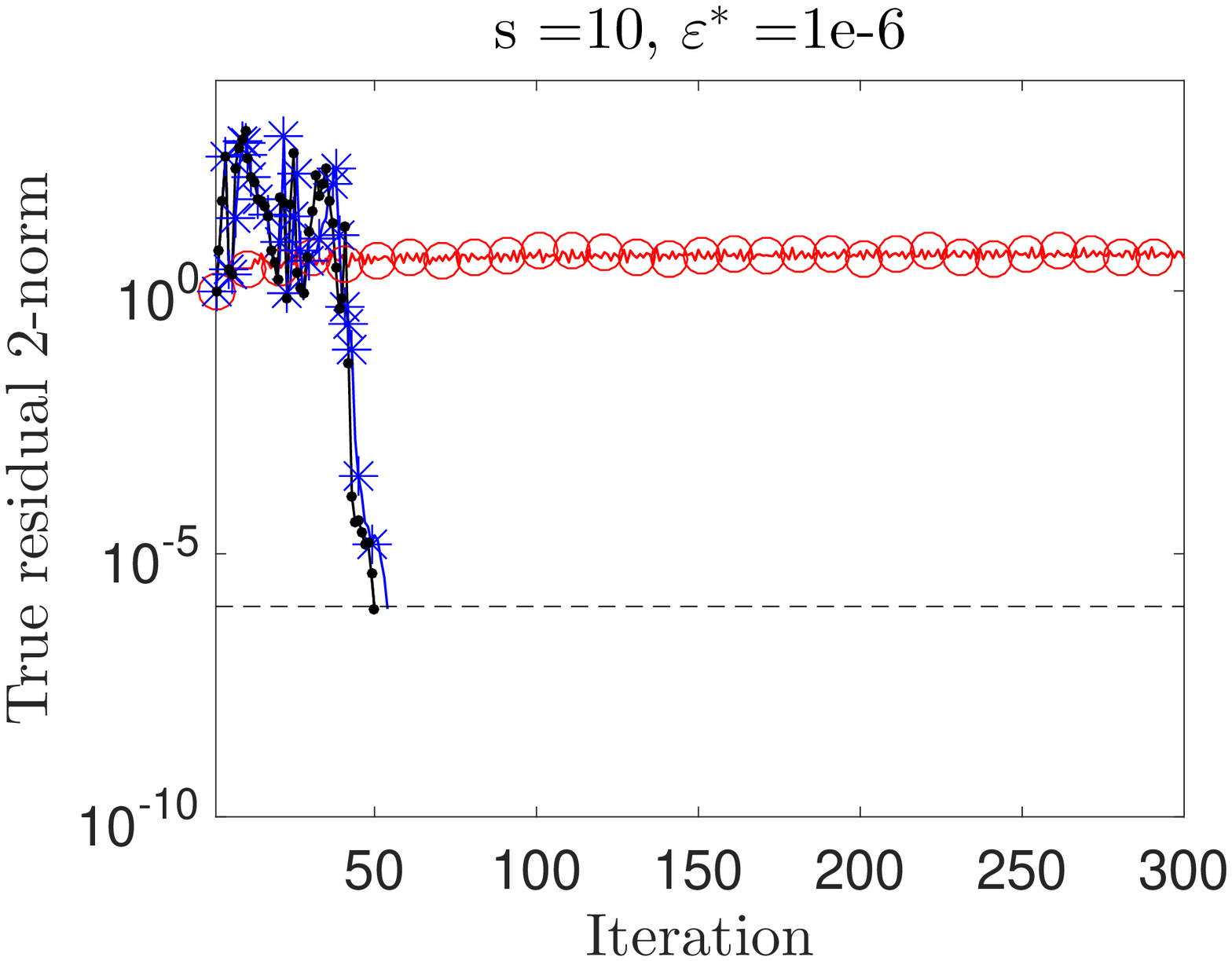}\\
	\label{fig:ex5}
	
	\vspace{.6cm}
	\footnotesize
\captionof{table}{Number of global synchronizations (outer loop iterations) in tests for matrix ex5 (corresponding to Figure~\ref{fig:ex5}). Dashes (-) in the table indicate that the method failed to converge to the desired level $\veps^*$. Tests used $c_{k,j}=\sqrt{\kappa(A)}$.}
\begin{tabular}{lr|c|c|c|}
\cline{3-5}
                                                      & \multicolumn{1}{l|}{} & fixed $s$-step CG & adaptive $s$-step CG & classical CG        \\ \hline
\multicolumn{1}{|l|}{\multirow{3}{*}{$\veps^*=$1e-8}} & $s=4$                 & 59          & 48               & \multirow{3}{*}{76} \\ \cline{2-4}
\multicolumn{1}{|l|}{}                                & $s=8$                 & -           & 45               &                     \\ \cline{2-4}
\multicolumn{1}{|l|}{}                                & $s=10$                & -           & 45               &                     \\ \hline
\multicolumn{1}{|l|}{\multirow{3}{*}{$\veps^*=$1e-6}} & $s=4$                 & 34          & 27               & \multirow{3}{*}{49} \\ \cline{2-4}
\multicolumn{1}{|l|}{}                                & $s=8$                 & -           & 27               &                     \\ \cline{2-4}
\multicolumn{1}{|l|}{}                                & $s=10$                & -           & 27               &                     \\ \hline
\end{tabular}
\label{tab:ex5}
\end{figure}

Our results are plotted in Figures~\ref{fig:gr3030}-\ref{fig:ex5}. 
For each test matrix, we ran classical CG (Algorithm~\ref{alg:cg}), as well as both (fixed) $s$-step CG (Algorithm~\ref{alg:scg}) and adaptive $s$-step CG (Algorithm~\ref{alg:vscg2}) for $s_{\text{max}}$ values $4$, $8$, and $10$ (for fixed $s$-step CG, $s=s_{\text{max}}$). For each matrix and $s_{\text{max}}$ value, we tested two different $\veps^*$ values: $\veps^*=\veps_{\text{CG}}$, the maximum attainable accuracy of classical CG, and $\veps^*=1e-6$. 
The monomial basis was used for both variable and fixed $s$-step methods in all tests (note that there is no technical restriction to the monomial basis; the analysis used in deriving the adaptive $s$-step CG method still applies regardless of the polynomials used in constructing $\yh_{k,\bar{s}_k}$). In all tests for adaptive $s$-step CG, we set $\bar{s}_0 = f = s_{\text{max}}$ (see Algorithm~\ref{alg:vscg2}), as this gives the optimal choice of $s_k$ in each iteration. (We also tried using smaller values of $f$, e.g., $f=1$ and $f=2$; this did not significantly change the total number of outer loop iterations performed). We used $c_{k,j}=1$ in our condition for determining the number of inner loop iterations to perform in all experiments except for tests with matrices bcsstk09 (Figure~\ref{fig:bcsstk09}) and ex5 (Figure~\ref{fig:ex5}); for these problems, convergence is more irregular, and thus the bound~\eqref{yhk} is tighter.  

Each figure has a corresponding table (Tables~\ref{tab:gr3030}-\ref{tab:ex5}) which lists the number of outer loop iterations (a proxy for the latency cost or number of global synchronizations) required for convergence of the true residual $2$-norm to level $\veps^*$. For classical CG, the number reported gives the total number of iterations. Dashes in the tables indicate that the true residual diverged or stagnated before reaching level $\veps^*$.

\begin{table}[ht]
\centering
\footnotesize
\caption{Test Matrix Properties.}
\label{mats}
\begin{tabular}{|l|l|l|l|l|}
\hline
Matrix     & $n$  & nnz   & $\norm{A}$ & $\kappa(A)$ \\ \hline
gr\_30\_30 & 900  & 7744  & 1.5        & 1.9e2       \\ \hline
mesh3e1    & 289  & 1377  & 1.8        & 8.6e0        \\ \hline
nos6       & 675  & 3255  & 2.0        & 3.5e6       \\ \hline
bcsstk09   & 1083 & 18437 & 2.0        & 1.0e4       \\ \hline
ex5        & 27   & 279   & 3.8        & 5.7e7       \\ \hline
\end{tabular}
\end{table}


\subsection{Observations}

In all experiments, the adaptive $s$-step CG method was able to converge to the desired tolerance $\veps^*$, whereas for $\veps^*=\veps_{\text{CG}}$ and higher $s$ values ($s=8$ and/or $s=10$), fixed $s$-step CG may not achieve convergence to this level. We point out that even when the matrix is very well-conditioned, choosing $s$ too large can very negatively affect the convergence rate in fixed $s$-step CG; see, e.g., the case with $s=10$ and $\veps^*=1e-6$ in Figure~\ref{fig:mesh3e1}, where $s$-step CG requires more than 8 times the number of synchronizations than classical CG. 
In tests where both fixed $s$-step CG and adaptive $s$-step CG converge, adaptive $s$-step CG takes at most $2$ additional outer iterations versus fixed $s$-step CG (although in most cases it takes as many or fewer outer iterations). We note that in all cases, variable $s$-step CG uses fewer outer loop iterations (synchronization points) than classical CG, so we still expect a performance advantage on latency-bound problems. 

The experiments using $\veps^*=1e{-}6$ demonstrate the ability of the adaptive $s$-step CG method to automatically adjust to use the largest $s_k$ values possible (according to our bound) that will eventually attain the desired accuracy. For example, in Figure~\ref{fig:gr3030} for gr\_30\_30, for all tests using $\veps^*=1e{-}6$, adaptive $s$-step CG uses $s_k=s_{\text{max}}$ in every outer iteration, defaulting to the fixed $s$-step method. In contrast, in plots on the left using $\veps^*=\veps_{\text{CG}}$, adaptive $s$-step CG uses smaller $s_k$ values at the beginning and then increases up to $s_k=s_{\text{max}}$ and is able to attain accuracy $\veps_{\text{CG}}$ whereas fixed $s$-step CG fails when for $s=8$ and $s=10$.

Below we list the sequence of $s_k$ values chosen in adaptive $s$-step CG for tests with $\veps^*=\veps_{\text{CG}}$ and $s_{\text{max}}=10$ (plots in the lower left in each figure):
\begin{itemize}
\item gr\_30\_30: 1, 1, 2, 2, 2, 3, 3, 3, 4, 5, 6, 8, 10, 10
\item mesh3e1: 1, 1, 2, 4, 6, 9, 10
\item nos6: 6, 1, 2, 3, 4, 5, 4, 4, 1, 4, 4, 5, 6, 7, 8, 10, 10, 10, 10, 10, 10, 10, 10, 10, 10, 10, 10, 10, 10, 10, 10, 10, 10, 10, 10
\item bcsstk09: 2, 2, 3, 3, 3, 3, 3, 1, 2, 3, 3, 3, 2, 2, 2, 2, 3, 3, 3, 3, 4, 4, 4, 4, 5, 5, 6, 7, 1, 6, 6, 7, 7, 7, 7, 8, 8, 8, 9, 9, 10, 10, 10, 10, 10, 10, 10, 10, 10, 10
\item ex5: 1, 1, 1, 1, 1, 1, 1, 1, 1, 1, 1, 1, 1, 1, 1, 1, 1, 3, 1, 1, 2, 1, 1, 1, 3, 3, 2, 1, 1, 1, 2, 1, 1, 2, 4, 6, 2, 4, 5, 7, 2, 6, 1, 7, 6
\end{itemize}
There are a few things to notice. First, as expected, the sequence of $s_k$ values is monotonically increasing when the residual 2-norm decreases monotonically or close to it (matrices gr\_30\_30 and mesh3e1), but when the residual 2-norm oscillates, we no longer have monotonicity of the $s_k$ values (matrices nos6, bcsstk09, and ex5). Second, it may be that $s_{\text{max}}$ is never reached, as is the case for matrix ex5. It is also worth pointing out that the sequence of $s_k$ values need not start with $s_0=1$; for nos6, $s_0=6$ and for bcsstk09, $s_0=2$.

It may be that on a particular machine for a particular matrix structure, it is not beneficial in terms of minimizing time per iteration to go above a small value of $s_{\text{max}}$, say $s_\text{max}=4$. In this case, as long as the matrix is reasonably well-conditioned and the accuracy required is not too high, the fixed $s$-step CG method can probably be used without too much trouble. However, even in this case the adaptive $s$-step CG method can provide reliability while still maintaining the reduction in synchronizations of the fixed $s$-step method, at the cost of potentially wasting a few matrix-vector multiplications in each outer iteration.  

Looking at the number of outer loop iterations in Tables~\ref{tab:gr3030}-\ref{tab:ex5}, we see that in our tests for adaptive $s$-step CG it is often the case that increasing the $s_{\text{max}}$ value past a certain point gives diminishing returns in terms of reducing synchronization. This may represent a fundamental limit of the possible performance using the monomial basis for these problems. In these cases it may be that better performance (higher $s_k$ values) can be achieved only by using a better-conditioned polynomial basis for constructing $\yh_{k,\bar{s}_k}$ (e.g., Newton or Chebyshev bases).

\section{Conclusions and future work}
\label{sec:conclusion}

In this work, we developed the adaptive $s$-step CG method, which is a variable $s$-step CG method where 
the parameter $s_k$ is determined automatically in each outer loop such that a user-specified accuracy is attainable. 
The method for determining $s_k$ is based on a bound on the residual gap within outer loop $k$, from which we derive a constraint on the condition number of the $s$-step basis matrix generated in outer loop $k$. The computations required for determining $s_k$ can be performed without introducing any additional synchronizations in each outer loop. 
Our numerical experiments demonstrate that the adaptive $s$-step CG method is able to attain up to the same accuracy as classical CG while still significantly reducing the number of outer loop iterations (i.e., synchronizations) performed. In contrast, the fixed $s$-step CG method can fail to converge to this accuracy when $s>4$. 

Additionally, the adaptive $s$-step CG method adjusts to use larger $s_k$ when an application requires a less accurate solution, automatically optimizing the tradeoff between speed and accuracy. In this way, a user can use the matrix structure and machine parameters to estimate the value $s_{\text{max}}$ that will minimize the time per iteration (e.g., by guessing, or offline auto-tuning), input the accuracy they require, and the adaptive $s$-step CG method will choose $s_k$ as large as it can (i.e., reducing as much synchronization as it can) while still converging to the desired level. This takes the burden of parameter selection off the user. Our adaptive $s$-step CG method thus makes significantly advances in improving reliability and usability in $s$-step Krylov subspace methods.

Our future work includes a high-performance implementation and performance experiments in order to compare classical CG, fixed $s$-step CG, and adaptive $s$-step CG on large-scale problems. Another potential improvement is to obtain a tighter bound to replace~\eqref{fb}, which could avoid the requirement of adjusting the function $c_{k,j}$ (which represents a multiplicative factor in the bounds~\eqref{fb} and~\eqref{yhk}) in Algorithm~\ref{alg:vscg2}. 

The adaptive $s$-step CG method is designed to maximize $s_k$ subject to a constraint on accuracy, but obviously it would also be useful to optimize $s_k$ to achieve the \emph{overall} lowest latency/synchronization cost. This would require a better understanding of the finite precision convergence behavior of both $s$-step and classical Krylov subspace methods, which is an active area of research.

We believe that this adaptive approach based on a bound on the residual gap can be extended to other $s$-step Krylov subspace methods as well. The same approach outlined here can be used in $s$-step BICG (see, e.g.,~\cite{carson2013avoiding}), which uses the same recurrences for the solution and updated residual. We conjecture that a similar approach applies to $s$-step variants of other recursively computed residual methods (see~\cite{greenbaum1995estimating}), like $s$-step BICGSTAB~\cite{carson2013avoiding}, as well as methods like $s$-step GMRES. In the case of $s$-step GMRES, it may be possible to simplify the adaptive algorithm: since the residual norm monotonically decreases, the sequence of $s_k$ values should monotonically increase. Such a sequence of $s_k$ values was used in the variable $s$-step GMRES method of Imberti and Erhel~\cite{imberti2016vary}.   
 
As a final comment, we note that in our analysis in section~\ref{sec:varscg}, we used equation~
\eqref{fb} to determine a bound on $\kappa(\yh_{k,s_k})$ that guarantees 
accuracy on the level $\veps^*$. We could have just as easily used equation~
\eqref{fb} to develop a mixed precision $s$-step CG method where $s$ is 
fixed, but the precision used in computing in outer loop $k$ varies. In other words, rearranging the inequality~\eqref{fb} we see that, given that 
we've computed an $s$-step basis with condition number $\kappa(\yh_{k,s})$, if we 
want to achieve accuracy $\veps^*$, we must have that
\[
\veps_k \leq \frac{\veps^*}{c_{k,j} \kappa(\yh_{k,s}) \max_{\ell\in\{0,\ldots,j\}} \norm{\rhat_{m+\ell+1}}}.
\]
Investigating the potential for mixed precision $s$-step Krylov subspace 
methods based on this analysis remains future work.


\end{document}